\documentclass[a4paper,11pt]{article}
\RequirePackage[T1]{fontenc}
\RequirePackage{fix-cm}
\usepackage{amsmath,amsfonts,amsthm}
\usepackage{graphicx} 
\usepackage[unicode, final, colorlinks=true]{hyperref} 
\usepackage{cite}
\usepackage{caption}
\usepackage{mathtools}
\usepackage{newtxtext}
\usepackage[varvw]{newtxmath}
\usepackage{bm}
\usepackage[T1]{fontenc}
\usepackage[utf8]{inputenc}
\usepackage[dvipsnames]{xcolor}
\usepackage{tikz, pgfplots}
\usepackage{tikzscale}
\usepackage{booktabs}
\usepackage{multirow}
\usepackage{array}
\usepackage{csquotes}
\usepackage{siunitx}
\usepackage{accents}

\pgfplotsset{compat=1.18}
\usepgfplotslibrary{groupplots}
\usepgfplotslibrary{fillbetween}
\usetikzlibrary{backgrounds}
\usetikzlibrary{arrows.meta}
\pgfplotsset{plot coordinates/math parser=false}
\usetikzlibrary{external}
\usepgfplotslibrary{patchplots}
\usetikzlibrary{arrows.meta,bending,chains,decorations.markings}
\providecommand{\keywords}[1]{\textit{Keywords:} #1}
\providecommand{\msc}[1]{\textit{2010 MSC:} #1}

\newtheorem{definition}{Definition}

\newtheorem{proposition}{Proposition}

\allowdisplaybreaks[4]
\newcommand{\R}{\mathbb{R}}

\newcommand{\ddt}{\partial_t}
\newcommand{\ddx}{\partial_x}

\newcommand{\kmh}{\unit[per-mode = symbol]{\kilo\metre\per\hour}}
\newcommand{\pkm}{\unit[per-mode = symbol]{\kilo\metre^{-1}}}

\begin{document}
\title{Data-Driven Models for Traffic Flow at Junctions}
\author{
  Michael~Herty$^{1}$ \and
  Niklas~Kolbe$^{1}$\footnote{Corresponding author}
}

\date{
  \small
  $^1$Institute of Geometry and Practical Mathematics,\\ RWTH Aachen University, Templergraben 55,\\ 52062 Aachen, Germany\\
   \smallskip
   {\tt \{herty,kolbe\}@igpm.rwth-aachen.de} \\
   \smallskip
   \today
}
\maketitle
%
%
%
%
%
%
%
%
%
%
%
%
\begin{abstract}
The simulation of traffic flow on networks requires knowledge on the behavior across traffic intersections. For macroscopic models based on hyperbolic conservation laws there exist nowadays many ad-hoc models describing this behavior. Based on real-world car trajectory data we propose a new class of data-driven models with the requirements of being consistent to networked hyperbolic traffic flow models. To this end the new models combine artificial neural networks with a parametrization of the solution space to the half-Riemann problem at the junction. A method for deriving density and flux corresponding to the traffic close to the junction for data-driven models is presented. The models parameter are fitted to obtain suitable boundary conditions for macroscopic first and second-order traffic flow models. The prediction of various models are compared considering also existing coupling rules at the junction. Numerical results imposing the data-fitted coupling models on a traffic network are presented exhibiting accurate predictions of the new models.
\par
\vspace{0.5em}
\noindent
\msc{35L65, 90B10, 90B20}
\par
\vspace{0.5em}
\noindent
\end{abstract}

\section{Introduction}\label{sec:intro}
Mathematical modeling of vehicular traffic flow is considered on diffent scales and we refer to~\cite{MR2834083, MR3875436, MR3220784} for recent reviews. Here, we are interested in continuum models that evolve aggregated quantities like the vehicle density or the mean vehicle velocity in space and time. Those models are also called fluid-dynamic or macroscopic models and a rich literature exists today (e.g., \cite{aw2000resur, MR2784693, MR2366138, Dag1995, DAGANZO2006396, garavello2006traffflownetwor, MR2239057, PhysRevE.48.R2335, PhysRevE.50.54, Lebacque1993LesMM, lighthill1955kinemwavesii, payne71, payne1979freflo, richards1956shockwaveshighw, underwood61, Zha2002}). Among the (inviscid) macroscopic models one typically distinguishes between first-order models given by scalar hyperbolic conservation laws and second-order models comprised of systems of strictly hyperbolic equations. The most prominent first-order model has been introduced by Lighthill and Whitham \cite{lighthill1955kinemwavesii} and Richards \cite{richards1956shockwaveshighw} and it is based on a fixed relation between mean car velocity and density. This has later been extended e.g., towards a general family of functions by Aw and Rascle \cite{aw2000resur} and Zhang \cite{Zha2002} leading to a second-order model.  

 In order to extend those macroscopic models road networks, coupling or boundary conditions are required. The appropriate junction modeling strongly influences the dynamics on the roads and has been the focus of recent research in the field, see e.g.,~\cite{BreCanGar2014,GarHanPic2016, garavello2006traffflownetwor}. It is clear, that the conservation of mass at the junction needs to hold, but from a mathematical point of view this condition is not sufficient to obtain a well-posed initial boundary value problem. Therefore, further conditions have been imposed to obtain unique coupling conditions. As an example, the non-negative flux at the junction might be distributed according to given ratios modeling the preferences of the drivers and the total flux through the junction maximized according to specified ratios \cite{CocGarPic2005}. Many studies focus on suitable coupling conditions for the LWR model, see \cite{CocGarPic2005,GarHanPic2016,garavello2006traffflownetwor,GoaGoeKol2016, HerKla2003, holden1995} while for second order ARZ-type models similar considerations have been conducted. Coupling conditions for ARZ models have been proposed e.g. in~\cite{BuliXing2020,GarHanPic2016, garavello2006traffawrascl, garavello2006traffflownetwor, HauBas2007, HerMouRas2006, HerRas2006, KolCosGoa2018, KolGoeGoa2017, LebMamHaj2008, SieMauMou2009}, whereas in~\cite{ColGoaPic2010,GarMar2017} the authors proposed coupling conditions for ARZ phase transition models. We also refer to \cite{göttlich2021seconordertraff} for a comparison and novel conditions. Typically, well-posedness results are obtained using the notion of (nodal or half-) Riemann solvers to classify consistent boundary conditions. The precise definition of the Riemann solver depends on the system at hand and we refer to \cite{GarHanPic2016,LebMamHaj2008,HerRas2006} and below for more details. 

 Depending on the level of detail, some models have been employed and tested against data. Recently, it has been proposed that macroscopic models provide suitable flexibility to incorporate on-line traffic data and, in particular, fundamental diagram data~\cite{mmpdata, MR2642079, fan2014comparawrasclzhang}. While those models are now widely used only little work has been conducted on validation with traffic data \cite{mmpdata, Amin2008MobileCU, MR2509919, herty2018}. This situation is even more pronounced for traffic flow on road networks. Our manuscript contributes to the development of macroscopic coupling conditions using vehicle trajectory data. We will provide a framework for the development of general models at the junction that are consistent with the requirements given by the, also possibly data-fitted, macroscopic models on the road. In particular, a new class of machine learning models is proposed that combine trained artificial neural networks with a parametrization of the solution space of the half-Riemann problem at the junction to substitute for the nonlinear coupling conditions. We compare numerically the predictions of the novel models with existing classical models taken from~\cite{garavello2006traffflownetwor} and~\cite{göttlich2021seconordertraff}. We not only compare the coupling conditions but also conduct long-time predictions for upstream and downstream traffic simulated through a finite-volume method.

 The outline of the manuscript is as follows: in Section~\ref{sec:trafficflowmodels} we provide an overview on traffic flow models on networks, focusing on first order models, introduce a description of coupling models in terms of Riemann solvers and discuss admissible coupling fluxes in case of a 2-to-1 junction. In Section~\ref{sec:couplingmodels} we introduce various coupling models and describe them in the case of a 2-to-1 network. Along with established coupling rules from the literature we introduce our new machine-learning approach. Section~\ref{sec:data} considers the handling of vehicle trajectory data in the example of an on-ramp on a freeway. We discuss a normalization procedure with respect to the coupling delay and the estimation of road dependent fundamental diagrams. After elaborating on the parameter estimation of the coupling models in Section~\ref{sec:fitting} we validate these models on the road network in Section~\ref{sec:validation}.

 \section{Networked Traffic Flow Models}\label{sec:trafficflowmodels}
 In this section we give a summary on the mathematical modeling of traffic flow at a junction, introduce our notion of coupling models and show some properties on Riemann solvers and admissible conditions that will be used in the modeling in the following sections. The modeling of general junctions is considered in Section \ref{sec:trafficflowmodelsgeneral} and properties of 2-to-1 junctions that apply to the particular case of an on-ramp are discussed in Section~\ref{sec:trafficflowmodelsonramp}.

\subsection{Macroscopic models and Riemann solvers}\label{sec:trafficflowmodelsgeneral}
\begin{figure*}[t]
\centering
\vspace{.5cm}
\includegraphics{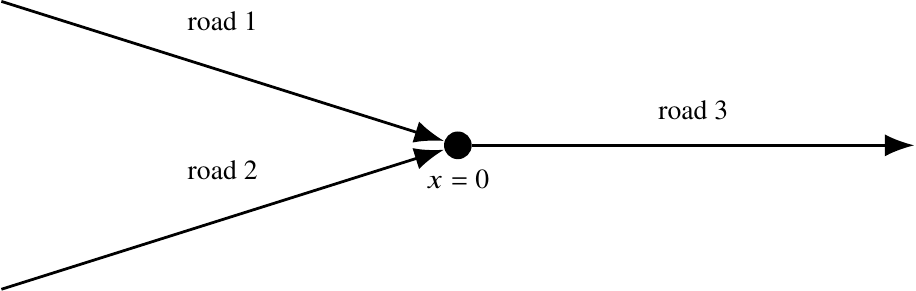}
\vspace{.5cm}
\caption{A road network consisting of three roads connected by a single junction.  We assume that road 1 denotes the entry lane, road 2 denotes a stretch of the freeway leading into the junction and road 3 denotes a stretch of the freeway going out of the junction. The roads are parameterized by the corresponding intervals $I_1=(-\infty,0]$, $I_2=(-\infty, 0]$ and $I_3=[0, \infty)$ assuming the junction at position $x=0$.}\label{fig:2to1}
\end{figure*}

Road networks are generally modeled by means of a directed graph with nodes representing the junctions and edges representing the roads, on which macroscopic traffic flow models, e.g.,~\eqref{eq:lwr} or~\eqref{eq:ARZmomentum} are imposed. While our approach can be generalized to road networks we focus here on the scenario of an on-ramp on a freeway modeled by a single junction with two incoming and one outgoing road. The three roads that meet at the 2-to-1 junction are the entry lane (road 1), the stretch of the freeway leading into the junction (road 2) and the stretch of the freeway after the merge (road 3) as shown in Figure~\ref{fig:2to1}. The junction itself is imposed at the physical location $x=0$, whereas the two incoming roads are parameterized by the intervals $I_1=I_2=(\infty, 0]$ and the outgoing one by the interval $I_3=[0, \infty)$.
Coupling conditions will be imposed at the junction that give rise to boundary conditions for the (system of) hyperbolic conservation laws. 

On each road $k\in\{1,2,3\}$ we describe traffic flow using the one dimensional Lighthill-Whitham-Richards (LWR) model~\cite{lighthill1955kinemwavesii, richards1956shockwaveshighw}, a hyperbolic conservation law that reads
\begin{equation}\label{eq:lwr}
\ddt \rho^k + \ddx \left( \rho^k \, V^k(\rho^k) \right) = 0 \quad \text{in }(0, \infty) \times I_k.
\end{equation}
 In this model $\rho^k=\rho^k(t,x)$ denotes the time and space dependent vehicle density on road $k$. The flux in \eqref{eq:lwr} is given by the density flow $f^k(\rho^k)=\rho^k \, V^k(\rho^k)$. A set of functions $\rho^k$ for $k \in \{1, 2, 3\}$ is called a \emph{weak solution} at the junction if for all families of test functions $\Phi^k:I_k \times (0,\infty) \to \mathbb{R}$ being smooth across the junction, i.e., $\Phi^1(t,0^-)=\Phi^2(t,0^-)= \Phi^3(t, 0^+)$ for all $t\geq0$, the following equation holds:
 \begin{equation} \label{eq:networkproblem1_LWR}
   \begin{split}
     \sum_{k \in \{1, 2, 3\}} \int_{0}^{\infty} \int_{I_k} \left[\rho^k  \partial_t \Phi^k + f^k \partial_x \Phi^k \right] \mathrm{d}x \, \mathrm{d }t 
     =  \sum_{k\in\{1,2, 3\}} \int_{I_k} \rho^{k,0} \Phi^k(0,x) \, \mathrm{d }x.
   \end{split}
 \end{equation}
In the above equation, $\rho^{k,0}$ denotes initial data. Note that the weak formulation \eqref{eq:networkproblem1_LWR} implies that conservation of mass is as necessary coupling condition, see e.g.,~\cite{holden1995}. The mass conservation condition is also called Kirchhoff condition and reads 
\begin{equation}\label{eq:kirchhofffullnet}
  \rho^1(t, 0^-) \, V^1(\rho^1(t, 0^-)) +  \rho^2(t, 0^-) \, V^2(\rho^2(t, 0^-)) =  \rho^3(t, 0^+) V^3(\rho^1(t, 0^+)), \qquad \text{for a.e.}\quad t > 0.
\end{equation}

In the LWR model, the velocity $V^k$ depends on both the road $k$ and the density $\rho^k$, the latter relation being known as fundamental diagram. We assume these diagrams to be chosen such that the flux functions $f^k$ for $k\in\{1, 2, 3\}$ are concave with respect to the density, satisfy $f^k(0)=f^k(\rho^k_\text{max})=0$ and have a unique global maximum $0<\sigma_k<\rho^k_\text{max}$ with $\rho^k_\text{max}$ denoting the maximal vehicle density on the road. A possible choice is the Greenshields model \cite{greenshields1935} 
\begin{equation}\label{eq:fundamentaldiagram}
  V^k(\rho^k) = v^k_\text{max} \left( 1- \frac{\rho^k}{\rho^k_\text{max}} \right),
\end{equation}
where the maximal/stagnation density $\rho^k_\text{max}$ and the maximal velocity $v^k_\text{max}$ are road dependent parameters. In the second order model by Aw, Rascle \cite{aw2000resur} and Zhang \cite{Zha2002} model \eqref{eq:lwr} is complemented by the additional equation 
\begin{equation}\label{eq:ARZmomentum}
  \ddt \left(\rho^k \, w^k\right) + \ddx\left( \rho^k \, w^k \,  V^k\right) = 0 \quad \text{in }(0, \infty) \times I_k,
\end{equation}
where $w^k$ denotes the Lagrangian marker on road $k$ carrying information about the traffic at a given point in space and time. The velocity in the ARZ model is determined by the marker and a pressure law $p(\rho^k, w^k)$ through the relation
\begin{equation}\label{eq:velocitysecondorder}
  V^k = w^k - p(\rho^k, w^k)
\end{equation}
replacing the fundamental diagram ~\eqref{eq:fundamentaldiagram} in the first order model. Choosing the pressure $p(\rho^k, w^k) = \frac{w^k} {\rho^k_\text{max}}\, \rho^k$ we recover relation~\eqref{eq:fundamentaldiagram} so that the variable marker $w^k$ takes the role of the maximal velocity $v^k_\text{max}$. Other choices for the pressure have been investigated e.g., in \cite{fan2014comparawrasclzhang}. A corresponding weak formulation for the momentum equation~\eqref{eq:ARZmomentum} states the conservation of the quantity $\rho w$ across the junction as necessary condition:
\begin{align*}
	\sum\limits_{i \in \delta^-_v}  \left( \rho^i w^i v^i \right)(t,b_i) = \sum\limits_{j  \in \delta^+_v}  \left( \rho^j w^j v^j \right)(t,a_j), \; t \geq 0. 
\end{align*}
 For both models, the relations accounting for conservation across the junction are not sufficient to obtain well-defined boundary conditions. 

In order to formalize the coupling condition in a unified presentation we will in the following use the notion of \emph{Riemann solvers} (RS), see e.g. \cite{garavello2006traffflownetwor}. Although those solvers are also available for the second-order models, see \cite{HerRas2006}, we will consider RS for the first-order models only. This is also due to the data available that does not yet include information on the Lagrangian marker $w^k$. However, motivated by the approximation qualities of the ARZ model, below we consider coupling conditions that also allow to modify the variable maximal velocity $v^k_\text{max}$ in \eqref{eq:fundamentaldiagram}. This corresponds to having a spatially constant marker $w^k$ whose value can be changed at the junction through the coupling condition. The details will be given in the forthcoming section.

\begin{definition}\label{def:rs}
  A \emph{Riemann solver} (RS) to the LWR system \eqref{eq:lwr} on the 2-to-1 network is a function
  \[
    \mathcal{RS}_\rho: \R^{3}_{\geq0} \rightarrow \R^{3}_{\geq0}
  \]
  that maps data $(\rho_0^{1-},\rho_0^{2-},\rho_0^{3+})$ accounting for initial conditions of a Riemann problem at the junction (being constant on each road) to the \emph{coupling data} $(\rho_R^1, \rho_R^{2},  \rho_L^{3})$. On the roads $I_1$, $I_2$ and $I_3$ the solution of the network system is given by the solution of the initial value boundary problem with initial data $\rho_0^i$ and Dirichlet boundary data $\rho_R^i$/$\rho_L^i$. Coupling data imposed by the RS is such that waves of the network solution have negative speed on incoming roads and positive speed on outgoing roads. The RS is \emph{consistent} if the following condition holds
  \[
    \mathcal{RS}_\rho \left(\mathcal{RS}_\rho(\rho_0^{1-},\rho_0^{2-},\rho_0^{3+})) \right) = \mathcal{RS}_\rho(\rho_0^{1-},\rho_0^{2-},\rho_0^{3+}).
    \]
\end{definition}

By its definition the RS provides boundary conditions for $\rho^k(t,0)$ under suitable assumptions. Then, using wave-front tracking, well-posedness of macroscopic first- and second-order models has been established, and we refer to~\cite{BreCanGar2014} for a recent survey on those techniques. In case of non-constant data on the roads the role of the Riemann data is transferred to traces of the densities at the coupling node. Under our assumptions on the flux functions $f^{k}$ the wave speed property of the RS is implied by the \emph{demand and supply conditions}~\cite{lebacque1996godun}
  \begin{equation}\label{eq:demandsupplyconditions}
    f_0^1 \leq d^1(\rho_0^{1-}), \qquad f_0^2 \leq d^2(\rho_0^{2-}),  \qquad   f^{3}_0  \leq  s^3(\rho_0^{3+}),
\end{equation}
where, denoting the density of maximal flow as $\sigma^k$, the demand and supply functions are defined by
 \begin{equation}\label{eq:demandsupply}
   d^k(\rho^k) =
   \begin{cases}
     f^k(\rho^k) &\text{if }\rho^k \leq \sigma^k ,\\
     f^k(\sigma^k) &\text{if }\rho^k > \sigma^k
   \end{cases} \quad \text{for }k\in \{1,2 \},\qquad
   s^3(\rho^3) =
   \begin{cases}
     f^3(\sigma^3) &\text{if }\rho^3 \leq \sigma^3, \\
     f^3(\rho^3) &\text{if }\rho^3 > \sigma^3.
   \end{cases}
 \end{equation}
 We note that the density flow obtained from the Greenshield model \eqref{eq:fundamentaldiagram} attains maximal flow at density $\sigma^k=\rho^k_\text{max}/2$.  Demand and supply have also been defined for ARZ-type models in~\cite{HerRas2006,LebMamHaj2008,göttlich2021seconordertraff}.

 As shown e.g., in~\cite{holden1995} the Kirchhoff condition and the notion of demand and supply are not sufficient to define unique boundary conditions for $\rho^k$ at the junction and hence the RS needs to incorporate additional conditions. To describe these in the following we use an alternative representation of a RS of the form
 \begin{equation}\label{eq:riemannsolver}
   \mathcal{RS}: \left(\rho_0^{1-},\rho_0^{2-},\rho_0^{3+}\right) \mapsto \left(f_0^1, f_0^2, f_0^3\right),
 \end{equation}
 which maps Riemann data to \emph{coupling fluxes} instead of coupling data in terms of densities as considered in Definition~\ref{def:rs}. This alternative form allows for a simplified and consistent description of the coupling models in Section \ref{sec:couplingmodels}. The output of \eqref{eq:riemannsolver} can be transformed to corresponding coupling data.
 
 \begin{proposition}\label{prop:couplingdata}
   Let $\mathcal{RS}: \R^{3}_{\geq0} \rightarrow \R^{3}_{\geq0}$ be a RS of the form \eqref{eq:riemannsolver} to the LWR system \eqref{eq:lwr} at the junction of the 2-to-1 network. Suppose that the flux functions $f^1$, $f^2$ and $f^3$ are concave with respect to the density, satisfy $f^k(0)=f^k(\rho^k_\text{max})=0$ and have a unique global maximum $0<\sigma_k<\rho^k_\text{max}$. Further assume that for given $\rho_0^{1-}$, $\rho_0^{2-}$ and $\rho_0^{3+}$ the coupling fluxes $f_0^{1}$, $f_0^{2}$ and $f_0^{3}$ satisfy the demand and supply conditions~\eqref{eq:demandsupplyconditions}. Then, defining the densities
   \[
     \rho^k_R =
     \begin{cases}
       \rho_0^k &\text{if }f_0^k=f^k(\rho_0^k) ,\\
       {( \overline{f^k})}^{-1}(f_0^k)& \text{otherwise}
     \end{cases}\quad \text{for }k\in \{1,2 \},\qquad
     \rho^3_L =
     \begin{cases}
       \rho_0^3 &\text{if }f_0^3=f^3(\rho_0^3) ,\\
       {(\underline{f^3})}^{-1}(f_0^3)& \text{otherwise}
     \end{cases}
   \]
   with $(\overline{f^k})^{-1}$ denoting the inverse function of $f^k$ on the interval $[\sigma^k, \rho^k_\text{max}]$ and $( \underline{f^k})^{-1}$ the inverse function of $f^k$ on the interval $[0, \sigma^k]$, a RS in the sense of Definition~\ref{def:rs} is obtained.  
 \end{proposition}
 \begin{proof} This result follows from the analysis in~\cite{holden1995}, the necessary details are provided here for convenience.
   A suitable RS is defined if for each edge $k\in\{1,2,3\}$ either $\rho_{R/L}^k= \rho_0^k$ or the Rankine-Hugoniot shock velocity at the junction,
   \[
     s^k_0 = \begin{cases}
               \frac{f_0^k - f^k(\rho_0^k)}{\rho_R^k - \rho_0^k} & k \in \{1,2\}, \\[5pt]
               \frac{f^k(\rho_0^k) -f_0^k}{\rho_0^k - \rho_L^k} & k =3
             \end{cases},
           \]
is nonpositive in incoming edges ($k=1$ or $k=2$) or nonnegative in the outgoing edge ($k=3$). We consider an incoming edge and exclude the trivial case $\rho_{R}^k= \rho_0^k$. If on the one hand $\rho_0^k>\sigma^k$ then by the demand condition in \eqref{eq:demandsupplyconditions} $\rho_R^k\in[\sigma^k, \rho^k_\text{max}]$ and due to $f_0^k=f^k(\rho_0^k)$ clearly $s^k_0$ is the secant slope of a strictly monotonically decreasing function and thus negative. If on the other hand $\rho_0^k<\sigma^k$ the demand condition requires $\rho_L^k\in ( {(\overline{f^k})}^{-1}(f^k(\rho_0^k)), \rho^k_\text{max}]$. Since $f^k$ is strictly monotonically decreasing in this interval this implies $f_0^k < f^k(\rho_0^k)$ and $\rho_R^k>\rho_0^k$ and therefore $s^k_0$ is also negative in this case. Nonnegativitity of $s^3_0$ is shown analogously exploiting the monotonicity of $f^3$ in the interval $[0, \sigma^3]$.
 \end{proof}
 
Consistency of an RS of the form \eqref{eq:riemannsolver} is defined analogously to Definition~\ref{def:rs} and is inherited when following the transformation in Proposition~\ref{prop:couplingdata}. Conversely, coupling data in terms of densities can be transformed to coupling fluxes by applying the flux functions corresponding to the respective roads. Existence and uniqueness of network solution for suitable coupling rules has been established in~\cite{holden1995, garavello2006traffflownetwor}. These results have also been extended to second-order models in e.g.,~\cite{HerRas2006,göttlich2021seconordertraff}.

\subsection{Admissible coupling fluxes}\label{sec:trafficflowmodelsonramp}
In this work we focus on coupling models of the form \eqref{eq:riemannsolver} determining coupling fluxes $f_0^{1}$, $f_0^{2}$ and $f_0^{3}$ that in addition to the demand and supply conditions~\eqref{eq:demandsupplyconditions} need to fulfill the mass conservation condition \eqref{eq:kirchhofffullnet} point-wise in time, i.e.,
    \begin{equation}\label{eq:kirchhoff}
  f_0^1 + f_0^2 = f_0^3.
\end{equation}
In the following, we analyze the set of admissible coupling fluxes.

\begin{definition}\label{def:G} For given trace data $\rho_0^{1-}$, $\rho_0^{2-}$ and $\rho_0^{3+}$ and concave flux functions $f^1$, $f^2$ and $f^3$ that satisfy $f^k(0)=f^k(\rho^k_\text{max})=0$ and have a unique global maximum $0<\sigma_k<\rho^k_\text{max}$ the admissible set $\mathcal{G} \subset \R^3$ refers to the set that includes all $(f_0^1, f_0^2, f_0^3)$ satisfying both the Kirchhoff condition \eqref{eq:kirchhoff} and the demand and supply conditions~\eqref{eq:demandsupplyconditions}. 
\end{definition}

\begin{figure}
  \centering
  \includegraphics{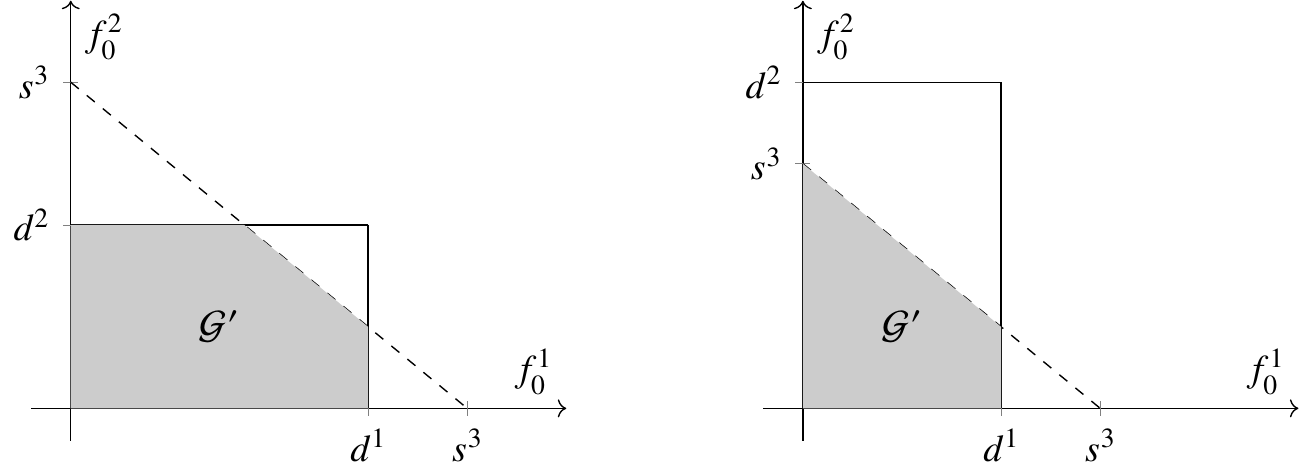}
  \caption{The restricted admissible set $\mathcal{G}^\prime$ in the two configurations $d^2<d^1<s^3<d^1+d^2$ (left) and $d^1<s^3<d^2$ (right).}\label{fig:G}
\end{figure}

Due to \eqref{eq:kirchhoff} any  $(f_0^1, f_0^2, f_0^3) \in \mathcal{G}$ has the form $(f_0^1, f_0^2, f_0^3)= (f_0^1, f_0^2, f_0^1 + f_0^2)$, thus we focus on the restricted set $\mathcal{G}^\prime=\{ (f_0^1, f_0^2):~(f_0^1, f_0^2, f_0^1 + f_0^2) \in \mathcal{G} \}$ to characterize the admissible set. Because of the demand and supply conditions~\eqref{eq:demandsupplyconditions} it holds
\begin{equation}\label{eq:gprime}
  \begin{split}
    (f_0^1, f_0^2) \in \mathcal{G}^\prime \qquad \text{iff} \qquad &0 \leq f_0^1 \leq \min\left\{ d^1(\rho_0^{1-}), s^3(\rho_0^{3+}) - f_0^2 \right\} \\ &\text{and} \quad 0 \leq f_0^2 \leq \min\left\{ d^2(\rho_0^{2-}), s^3(\rho_0^{3+}) - f_0^1 \right\},
  \end{split}
\end{equation}
hence, $\mathcal{G}^\prime$ is a rectangle in $\R^2$ that is (possibly) cut by a line as shown in Figure~\ref{fig:G} for two cases, see also \cite[Section 5.2.2]{garavello2006traffflownetwor}, where the flow maximization approach (see C1 below) is geometrically illustrated. We note that the admissible set depends only on the demand and supply at the junction and therefore write $\mathcal{G}=\mathcal{G}(d^1, d^2, s^3)$.  We state two properties of $\mathcal{G}$ in the following proposition.

\begin{proposition} The set $\mathcal{G}$ is non-empty and there is a bijective function $X_{\mathcal{G}}:[0,1]^2\rightarrow \mathcal{G}^\prime$.
\end{proposition}
\begin{proof}
  The first statement holds as $(0,0) \in \mathcal{G}^\prime$ for all trace data. A bijection is constructed as follows
  \begin{equation}\label{eq:X}
X_ {\mathcal{G}}: (\theta_1, \theta_2) \mapsto (\theta_1 \min\{d^1, s^3\}, \theta_2 \min\{d^2, s^3 - \theta_1 \min\{d^1, s^3\}\}),
\end{equation}
where we have neglected the arguments of the demand and supply functions for brevity.
\end{proof}

In the following we assume that the velocity on the three roads of the network is given by the Greenshields model \eqref{eq:fundamentaldiagram}. Motivated from \eqref{eq:velocitysecondorder} in the second order traffic model, we consider the parameterized flux functions
\begin{equation}\label{eq:parameterizedflux}
  f^k(\rho^k, w^k) = w^k \, \rho^k  \left( 1 - \frac{\rho^k}{\rho^k_\text{max}}\right), \qquad w^k \in [0, v^k_\text{max}]
\end{equation}
defined, such that $f^k(\rho^k, v^k_\text{max})=f^k(\rho^k)$. Including this parametrization by the Lagrangian marker also in \eqref{eq:demandsupply} gives rise to the generalized demand and supply functions $d^1(\rho^k, w^k)$, $d^2(\rho^k, w^k)$ and $s^3(\rho^k, w^k)$. The following result motivates the construction of our second order coupling models in the next section.

\begin{proposition}\label{prop:generalizedds}
  Let $w^k \in [0, v^k_\text{max}]$ for $k \in \{1,2,3\}$ and suppose that $\mathcal{RS}^*: \R^{3}_{\geq0} \rightarrow \R^{3}_{\geq0}$ is of the form \eqref{eq:riemannsolver} and such that the generalized demand and supply conditions
  \begin{equation}
    f_0^1 \leq d^1(\rho_0^{1-}, w^1), \qquad f_0^2 \leq d^2(\rho_0^{2-}, w^2),  \qquad   f^{3}_0  \leq  s^3(\rho_0^{3+}, w^3),
  \end{equation}
  as well as \eqref{eq:kirchhoff} hold. Then $\mathcal{RS}^*$ defines a RS.
\end{proposition}
\begin{proof}
  As the generalized demand and supply functions are monotonically nondecreasing in the marker parameter $w^k$ it holds
  \[
    \mathcal{G}(d^1(\rho_0^{1-}, w^1), d^2(\rho_0^{2-}, w^2), s^3(\rho_0^{3+}, w^3)) \subseteq  \mathcal{G}(d^1(\rho_0^{1-}), d^2(\rho_0^{2-}), s^3(\rho_0^{3+}))
  \]
  by the characterization \eqref{eq:gprime}, see also Figure~\ref{fig:G}. Therefore the original condition \eqref{eq:demandsupplyconditions} is implied so that $\mathcal{RS}^*$ satisfies the conditions in Proposition~\ref{prop:couplingdata} and thus defines a RS.
\end{proof}

\section{Coupling Models}\label{sec:couplingmodels}
In this section we recall various coupling models that impose coupling conditions of the type~\eqref{eq:riemannsolver} and satisfy both the Kirchhoff condition~\eqref{eq:kirchhoff} and the demand and supply conditions \eqref{eq:demandsupplyconditions} by construction. We distinguish between \emph{classical models} derived from literature conditions that are discussed in Section~\ref{sec:classic} and novel machine learning models that are discussed in Section~\ref{sec:mlmodels}. All models are also summarized in Table~\ref{tab:couplingmodels}. The models are stated for the 2-to-1 junction depicted in Figure~\ref{fig:2to1} for convenience. However, the presented framework can be generalized to an arbitrary junction. It is supposed that on all three roads the velocity is governed by the Greenshields model~\eqref{eq:fundamentaldiagram}.

\subsection{Classical Coupling Models}\label{sec:classic}
We consider four different classical coupling conditions also for comparison with the novel machine learning models later.
The models assume fixed preferences of the drivers with respect to their distribution at the junction, which is reflected by model parameters. For those as well as the Lagrangian markers used in our second order models parameter fitting is considered and described in Section~\ref{sec:fitting}.

\subsubsection*{C1: first order flow maximization}
A common approach discussed in \cite{garavello2006traffflownetwor} is the maximization of the flux in the junction such that the Kirchhoff condition \eqref{eq:kirchhoff} and the demand and supply conditions \eqref{eq:demandsupplyconditions} are fulfilled. Given the traces $\rho_0^{1-}$, $\rho_0^{2-}$ and $\rho_0^{3+}$  the demand in the incoming roads and the supply in the outgoing road according to \eqref{eq:demandsupply} is computed. Only, two cases can occur: In the case $ d^1(\rho_0^{1-}) +  d^2(\rho_0^{2-}) \leq s^3(\rho_0^{3+})$ the total flux is limited by the demand and maximal flux is obtained by taking
\begin{equation}\label{eq:lowdemand}
  f_0^1 =  d^1(\rho_0^{1-}), \qquad   f_0^2 = d^2(\rho_0^{2-}), \qquad f_0^3= d^1(\rho_0^{1-}) +  d^2(\rho_0^{2-}).
\end{equation}
In the contrary case $d^1(\rho_0^{1-}) +  d^2(\rho_0^{2-}) > s^3(\rho_0^{3+})$, the supply limits the total flux and an additional \emph{right of way} parameter $\beta\in[0,1]$ is required to determine the (preliminary) coupling fluxes according to
\begin{equation}\label{eq:distrule}
  \tilde f_0^1 = \beta \, s^3(\rho_0^{3+}), \qquad  \tilde f_0^2 = (1-\beta) \, s^3(\rho_0^{3+}), \qquad f_0^3 = s^3(\rho_0^{3+}).
\end{equation}
It might occur that one of the preliminary coupling fluxes $\tilde f_0^1$ and $\tilde f_0^2$ violates the demand condition. Suppose this happens in the first road, i.e., we have $\tilde f_0^1 > d^1(\rho_0^{1-})$, then we use $f_0^1=d^1(\rho_0^{1-})$ and $f_0^2= s^3(\rho_0^{3+}) - d^1(\rho_0^{1-})$. In the case $\tilde f_0^2 > d^2(\rho_0^{2-})$ the same procedure is applied with reversed indices 1 and 2. Otherwise the preliminary incoming coupling fluxes determine the final ones, i.e. $f_0^1=\tilde f_0^1$ and $f_0^2=\tilde f_0^2$. This procedure extends to the general case of multiple incoming and outgoing roads, see e.g.~\cite{CocGarPic2005}.

\subsubsection*{C2: second order flow maximization}
This model generalizes C1 and allows for changes in the fundamental diagrams. In more detail, we assume that in the junction the fluxes take the form~\eqref{eq:parameterizedflux} for some constant marker parameter $0 \leq w^k\leq v_\text{max}^k$. Note that due to Proposition~\ref{prop:generalizedds} this leads to a RS of the original coupled problem. Since $w^1$, $w^2$ and $w^3$ may be pre-computed during parameter estimation, the computation of the coupling fluxes follows the steps of C1 considering the updated fundamental diagrams. This procedure could also be used with spatially varying values of $w^k$ leading to conditions for the ARZ model as discussed in~\cite{HerRas2006, garavello2006traffawrascl}.

\subsubsection*{C3: alternative flow maximization}
While models C1 and C2 prioritize flow maximization over fixing a flux distribution as in \eqref{eq:distrule}, some publications such as \cite{göttlich2021seconordertraff, garavello2006traffawrascl} take the converse approach. Likewise, this model maximizes the fluxes subject to the constraints \eqref{eq:kirchhoff}, \eqref{eq:demandsupplyconditions} and \eqref{eq:distrule}. Additionally, this model employs the modified flux functions~\eqref{eq:parameterizedflux} in the junction and allows for changes of the parameter $w^k$ as considered in model C2. Its coupling fluxes are consequently given by
\begin{equation}
   f^1_0= \beta\, f_0^3, \quad f^2_0= (1-\beta)\, f_0^3 \quad f_0^3 = \min\left\{ \frac{d^1(\rho_0^{1-}, w^1)}{\beta}, \frac{d^2(\rho_0^{2-}, w^2)}{1-\beta},  s^3(\rho_0^{3+}, w^3)\right\}
\end{equation}
for an priority parameter $\beta\in(0,1)$. Due to the prioritized flux distribution all coupling fluxes become zero whenever $d^1=0$, $d^2=0$ or $s^3=0$, which is not the case in the maximization approach in models C1 and C2.

\subsubsection*{C4: homogenized pressure model}
In~\cite{göttlich2021seconordertraff} another coupling approach has been considered. If the coupling of model C3 is applied to the ARZ model, then, the priority will influence also the pressure $p(\rho,w)$ on the outgoing road. The computation of the homogenized pressure $p^\dag$ composed of the pressures corresponding to the mixture of $w^1$ and $w^2$, requires in general, to solve a nonlinear problem. Therefore, in \cite{göttlich2021seconordertraff} an approach based on the linearization of $p^\dag$ has been presented, that shows to be sufficiently close to the homogenized pressure and that governs the flux computation of this model. The coupling between the fluxes is as in model C3 but the dynamics introduced by modifying $p$ may lead to different boundary conditions for the densities $\rho^k$. 

\vspace{.5cm}
\begin{proposition} The coupling models C1, C2, C3 and C4 each define a consistent RS of the form~\eqref{eq:riemannsolver} that satisfies the demand and supply conditions~\eqref{eq:demandsupplyconditions} and therefore gives rise to suitable coupling data at the node using the construction in Proposition~\ref{prop:couplingdata}.
\end{proposition}

\subsection{Machine Learning Coupling Models}\label{sec:mlmodels}
In this section we introduce new machine learning models for the coupling of the network problem from Section~\ref{sec:trafficflowmodels}, which aim to mirror dynamics observed in real-word. To this end a significantly larger set of parameters is required in comparison to the models introduced in Section~\ref{sec:classic}. Besides the model architecture, on which we elaborate in this section, the \emph{training} of the parameters plays a crucial role in the model performance, see Section~\ref{sec:fitting} for details.

Neural network based machine learning models consists of multiple mappings (the \emph{layers}) that are applied successively and in parallel to assign a given model input to a suitable model output. In part, the layers are highly parameterized to allow the models to be fit to data in the training. Our models are constructed such that the predicted fluxes are included in the set of admissible coupling fluxes defined in Definition~\ref{def:G}. To achieve this, various problem specific layers are designed using the theoretical results in Section~\ref{sec:trafficflowmodels} and combined with an artificial neural network (ANN). In the following the construction is described for the 2-to-1 on-ramp junction but generalizations are straightforward.

\begin{figure}
  \centering
  \includegraphics[scale=0.88]{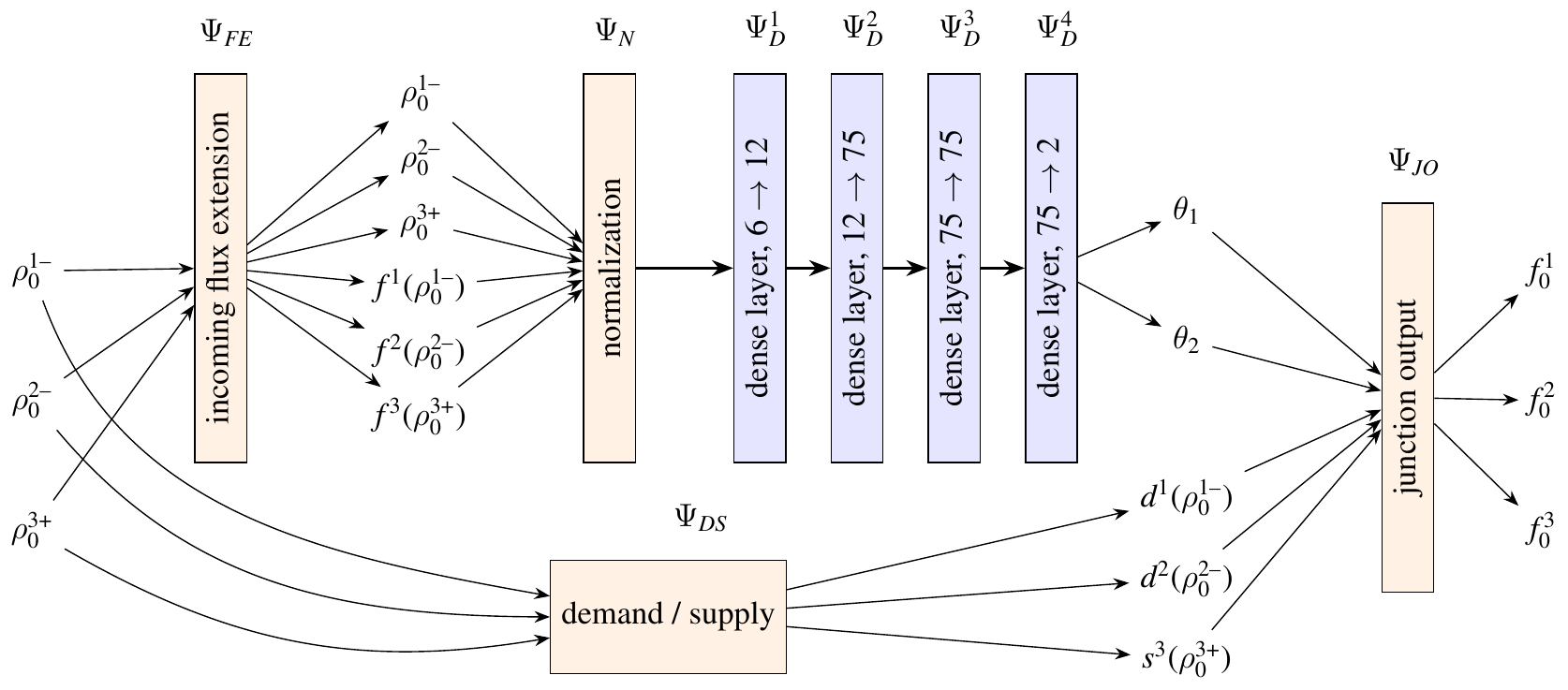}
\caption{Machine learning model architecture for ML2 and ML3. Trace data are first extended by incoming fluxes then normalized and processed by a four layer ANN determining the two parameters $\theta_1$ and $\theta_2$. In parallel demand and supply are computed from the model input, which together with the parameters are used to determine the output fluxes in the last layer of the model. In ML1 the trainable ANN part of the model (blue) is replaced by a single dense layer.}\label{fig:ML}
\end{figure}

\subsubsection{Problem specific layers}
Let the trace data  $\rho_0^{1-}$, $\rho_0^{2-}$, $\rho_0^{3+}\geq 0$ be given. Then we define the \emph{flux extension layer} as the map adding the incoming fluxes at the junction to the input data, i.e.,
\begin{equation}
  \psi_{FE}: (\rho_0^{1-}, \rho_0^{2-}, \rho_0^{3+}) \mapsto (\rho_0^{1-}, \rho_0^{2-}, \rho_0^{3+}, f^1(\rho_0^{1-}), f^2(\rho_0^{2-}), f^3(\rho_0^{3+})).
\end{equation}
Similarly, the \emph{demand and supply layer} evaluates the functions \eqref{eq:demandsupply} for given vehicle densities at the junction, i.e., 
\begin{equation}\label{eq:dslayer}
  \psi_{DS}: (\rho_0^{1-}, \rho_0^{2-}, \rho_0^{3+}) \mapsto (d^1(\rho_0^{1-}), d^2(\rho_0^{2-}), s^3(\rho_0^{3+})).
\end{equation}
For the model output we first determine the admissible set $\mathcal{G}$ from the demand and supply at the junction and then select a state within this set with the help of an ANN. To this end the mapping \eqref{eq:X} serves as a parametrization of the set $\mathcal{G}^\prime$. For given demand and supply obtained by the layer \eqref{eq:dslayer} and parameters $0\leq \theta_1, \theta_2 \leq 1$ goverened by the ANN, see below, we define the \emph{junction output layer}
\begin{equation}\label{eq:outputlayer}
  \Psi_{JO}: [0,1]^2 \times \R_{\geq 0}^3\times \rightarrow \R_{\geq0}^3, \qquad (\theta_1, \theta_2, d^1, d^2, s^3) \mapsto (f_0^1, f_0^2, f_0^3)
\end{equation}
such that $(f_0^1, f_0^2) = X_{\mathcal{G}}(\theta_1, \theta_2)$ for $\mathcal{G}=\mathcal{G}(d^1, d^2, s^3)$ and $f_0^3 = f_0^1 + f_0^2$. The problem specific layers are parameter free, i.e., they cannot be trained and its output depends only on the input obtained from previous layers of the model.

\subsubsection{Artificial neural network layers}
Dense ANN layers successively apply a linear affine map and a component-wise non-linear function (\emph{activation function})~\cite{pinkus1999approxmlp, marwala2019handb}. We employ such layers that map $m_j$ input neurons to $n_j$ output neurons and have the form
\begin{equation}\label{eq:denselayer}
  \Psi^j_D: \R^{m_j} \rightarrow [0, 1]^{n_{j}}, \qquad x \mapsto \sigma\left( W^j x + b^j \right),
\end{equation}
where $W^j\in\R^{m_j \times n_j}$ and $b^j \in \R^{n_j}$ are the weight and bias parameters that are fit in the training of the model. The sigmoid activation function that is used in \eqref{eq:denselayer} is component-wise given by
\begin{equation*}
\sigma(x) = \frac{1}{1 + e^{-x}}.
\end{equation*}
We consider fully connected ANNs that constitute compositions of $N$ such dense layers, i.e.,
\begin{equation}
\Psi_{ANN} = \Psi_D^N \circ \Psi_D^{N-1} \circ \dots \Psi_D^1
\end{equation}
such that $n_j=m_{j+1}$ for $j=1,\dots,N-1$. Our model architecture moreover requires $m_1=6$ and $n_N=2$.

\subsubsection{Model architecture}
In our machine learning approach we combine the defined layers to obtain coupling models of the form \eqref{eq:riemannsolver}. In parallel we 1) add incoming fluxes to the trace data, normalize the resulting vector and pass it to an ANN and 2) compute demand and supply from the trace data. The results are eventually used to compute coupling data using the mapping~\eqref{eq:outputlayer}. We refer to Figure~\ref{fig:ML} for a diagram of the procedure in the particular case of ML2 and ML3, see below for their definition. Formally our machine learning models are thus mappings of the form
\begin{equation}\label{eq:ml}
\mathcal{ML}: (\rho_0^{1-}, \rho_0^{2-}, \rho_0^{3+}) \mapsto \Psi_{JO}(\Psi_{ANN} \circ \Psi_N \circ \Psi_{FE}(\rho_0^{1-}, \rho_0^{2-}, \rho_0^{3+}), \Psi_{DS}(\rho_0^{1-}, \rho_0^{2-}, \rho_0^{3+})).
\end{equation}
Here $\Psi_N: \R^6 \to \R^6$ refers to the normalization layer, which is component-wise affine linear using parameters computed from the training data that are not trainable, see Section \ref{sec:fitting} for further details. The three model variants that we consider differ in the choice of the ANN architecture and their training.

\subsubsection*{ML1: single layer machine learning model}
To study the complexity of the problem we consider at first a machine learning model of the form \eqref{eq:ml}, in which $\Psi_{ANN}$ consists only of a single layer. In other words, the ANN part of the model here has the form
\begin{equation}
  \Psi_{ANN}(x) = \sigma\left(W x + b \right), \qquad W\in\R^{2 \times 6}, \qquad b \in \R^2.
\end{equation}
This model has 14 parameters to be determined in the training. To this end a loss function computing the mean squared error with respect to the data is minimized over the parameter space, i.e., a minimizer of the  optimization problem 
\begin{equation}\label{eq:noconstraints}
    \begin{aligned}
      &\min_{\mathcal{P}} \qquad & \|(f_d^1, f_d^2, f_d^3) - (f_0^1, f_0^2, f_0^3) \|_{I^h_\text{train}}^2 
    \end{aligned}
  \end{equation}
is sought, where the index $d$ indicates data and $\mathcal{P}$ the model parameters, being the entries of the matrix $W$ and the vector $b$ in this case. For details on the data and the used norm we refer to Section~\ref{sec:data}. The numerical solution of \eqref{eq:noconstraints} is explained in Section~\ref{sec:fitting}.

\subsubsection*{ML2: four layer machine learning model}
This variant of our machine learning approach increases the model complexity compared to ML1. This is achieved by employing the four layer ANN
\begin{equation}
  \begin{split}
    \Psi_{ANN} = \Psi_D^4 \circ \Psi_D^3 \circ \Psi_D^2 \circ \Psi_D^1, \qquad \Psi_D^1:\R^6 \to \R^{12}, \quad \Psi_D^{2}: \R^{12} \to \R^{75},\\
    \Psi_D^{3}: \R^{75} \to \R^{75}, \quad \Psi_D^{4}: \R^{75} \to \R^{2},
  \end{split}
\end{equation}
where each of the dense layer is of the form \eqref{eq:denselayer}.
The number of layers, input and output neurons has been chosen so that good training performance was obtained when testing the model capability, see Section~\ref{sec:fitting}. In this case the ANN has 6.911 parameters 
that are trained solving the optimization problem~\eqref{eq:noconstraints}. Figure~\ref{fig:ML} shows a diagram of the machine learning model~\eqref{eq:ml} in this variant.

\subsubsection*{ML3: consistent four layer machine learning model}
Lastly, we introduce a modification of the four layer model ML2 with a modification of the training. The approach~\eqref{eq:ml} guarantees a coupling model/RS that satisfies both conditions \eqref{eq:demandsupplyconditions} and \eqref{eq:kirchhoff}. Unlike the classical models introduced in Section~\ref{sec:classic} the machine learning models are not necessarily consistent RS in the sense of Definition~\ref{def:rs}. This model variant aims for consistency by including it as a training constraint. Given the traces $\rho_d^1$, $\rho_d^2$ and $\rho_d^3$ together with the corresponding fluxes $f_d^1$, $f_d^2$ and $f_d^3$ in the training data we consider the machine learning model output $(f_0^1, f_0^2, f_0^3) \coloneq \mathcal{ML}( \rho_d^1, \rho_d^2, \rho_d^3)$. Corresponding to the model output we compute the coupling data $(\rho_0^1, \rho_0^2, \rho_0^3)$ following Proposition~\ref{prop:couplingdata}. The four layer model is then trained solving the constrained problem
\begin{equation}\label{eq:consistencyconstraints}
    \begin{aligned}
      &\min_{\mathcal{P}} \qquad & \|(f_d^1, f_d^2, f_d^3) - (f_0^1, f_0^2, f_0^3) \|_{I^h_\text{train}}^2 \\
       &\text{subject to } \qquad & (f_0^1, f_0^2, f_0^3) = \mathcal{ML}( \rho_0^1, \rho_0^2, \rho_0^3).
    \end{aligned}
  \end{equation}
Details about the numerical solution of \eqref{eq:consistencyconstraints} are given in Section~\ref{sec:fitting}. We note that the additional constraint in the training only leads to an approximate consistency in the resulting machine learning model.

\begin{table}\fontsize{7}{7}\selectfont
  \centering
  \caption{Summary of the considered coupling models. Table entries indicate the model parameters and whether the consistency condition is satisfied.}\label{tab:couplingmodels}
  \begin{tabular}{c l c r l l}
    \toprule
    \multicolumn{2}{c}{model} & consistency & \multicolumn{2}{c}{parameters} & remark\\ \midrule
    C1 & first order flow maximization &  yes & 1: & $\beta$ & proposed in \cite{garavello2006traffflownetwor}\\
    C2 & second order flow maximization &  yes & 4: & $\beta$, $v^k_\text{max}$ & generalization of C1\\
    C3 &  alternative flow maximization &  yes & 4: & $\beta$, $v^k_\text{max}$ & proposed in \cite{garavello2006traffawrascl}\\
    C4 &  homogenized pressure model &  yes & 4: & $\beta$, $v^k_\text{max}$ & proposed in \cite{göttlich2021seconordertraff}\\
    ML1 &  single layer machine learning model & no & 14: & $W$, $b$ & single layer variant of Figure~\ref{fig:ML}\\
    ML2 &  four layer machine learning model & no & 6911: & $W^j$, $b^j$ for $j\in\{ 1,2,3,4\}$ & model from Figure~\ref{fig:ML}\\
    ML3 &  consistent four layer machine learning model & approximately & 6911: & $W^j$, $b^j$ for $j\in\{ 1,2,3,4\}$ & ML2 with training constraints\\
    \bottomrule
  \end{tabular}
\end{table}
\section{Traffic Data}\label{sec:data}
\begin{figure}[t]
  \centering  
  \includegraphics{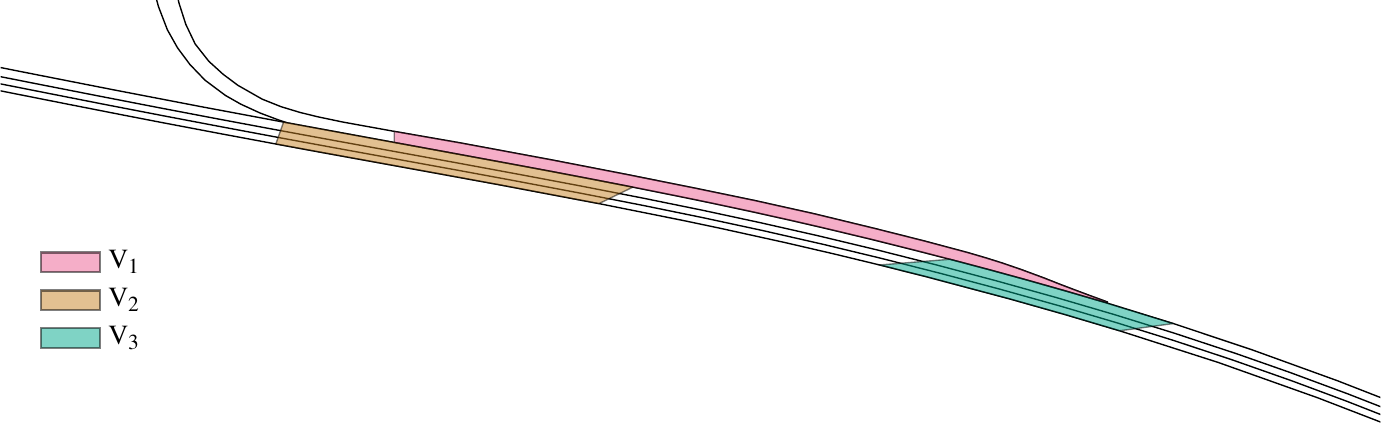}
  \caption{Bonn-Beuel entry lane with control volumes. Traffic either passes the three main lanes from left (east) to right (west) or enters from the on-ramp north. The east-west axis is inverted for consistency with Figure \ref{fig:2to1}.}\label{fig:laneBB}
\end{figure}

We consider car trajectory data, which provides a Lagrangian description of vehicle movement, at the German freeway A565 near the interchange Bonn-Beuel collected in May 2019. Traffic on and around an on-ramp was recorded by drone photography and vehicle trajectories were obtained from image processing of the recordings and then interpolation by third order polynomials. The recorded freeway stretch, shown in Figure~\ref{fig:laneBB}, covers approximately 270 meters and (including the on-ramp) 4 unidirectional lanes. In total, 31 data sets were recorded over 4 days either in the morning or in the afternoon, each one covering the traffic over periods of approximately 5 minutes in average. The number of vehicles entering the junction varied with the time of day, in average 44 vehicles passed the junction per minute throughout the data. More details on the datasets including the time of recording and average speeds are given in Table~\ref{tab:datasets}.\footnote{The data is available upon request.} Each data set $\mathcal{D}_k$ is represented by a family of smooth curves,  $\mathbf x_j: I_\ell \rightarrow \R^2$ with $j\in \mathcal{D}_\ell$ and $I_\ell\subset \R$ denoting the relevant time interval, describing the position of the vehicle barycentres with respect to time. We note that the time intervals $I_\ell$ are disjoint and denote the collection of all data by $\mathcal D = \cup_{\ell=1}^{31} \mathcal{D}_\ell$.

To obtain macroscopic coupling information we define specific regions/volumes at the junction as shown in Figure~\ref{fig:laneBB}. The first volume $V_1$ is passed by vehicles entering the junctions from the on-ramp, the second volume $V_2$ is passed by the transit traffic and the third volume $V_3$ is passed by the vehicles coming out of the junction. Traffic from the on-ramp can enter the freeway by changing to the northern main lane while passing a long stretch of the on-ramp, the position of lane change depends on various factors, such as vehicle speed and the current transit traffic. For this reason the volume $V_1$ covers a longer stretch of road than the volumes $V_2$ and $V_3$ and we assume a gap, i.e., a stretch of road that is included in neither volume, between the volumes $V_2$ and $V_3$ to guarantee that the volume $V_2$ is only passed by transit traffic.

From the data we can compute at each time the empirical vehicle density in volume $k$ by
\begin{equation}\label{eq:density}
   \rho_d^k(t) = \operatorname{card}\{ j \in \mathcal{D}:~\mathbf x_j(t)\in V_k  \} \operatorname{diam}(V_k)^{-1},
 \end{equation}
 where formally completed trajectories $\mathbf x_j: \R \rightarrow \R^2\cup\{ \infty\}$ such that $x_j(t) =\infty$ for all $t\notin I_\ell$ and $j\in \mathcal{D}_\ell$ are considered. These quantities serve as Riemann data input to the coupling models introduced in Section~\ref{sec:couplingmodels}. We define the time dependent empirical velocity in volume $k$ by computing the averages 
\begin{equation}\label{eq:velocity}
  v_d^k(t) =
  \begin{cases}
  \operatorname{card}{\{ j \in \mathcal D:~\mathbf x_j(t)\in V_k \}}^{-1} \sum_{\mathbf x_j(t) \in \mathcal D} \| {{\mathbf x}}_j^\prime (t)\|_2 &\text{if } \rho_d^k(t)>0, \\
    0, & \text{if } \rho_d^k(t) = 0.
  \end{cases}
  \end{equation}
  As the vehicle trajectories are given by polynomials, the derivative $\mathbf x_j^\prime(t)$ used in \eqref{eq:velocity} is well-defined for any $t \in  I_\ell$ assuming $j \in\mathcal{D}_\ell$. The empirical velocities allow us to define the empirical traffic flows $f_d^k =  \rho_d^{k}\,  v_d^k$, which are used to fit the coupling models. The introduced quantities will be considered on the discretized time intervals $I^h_\ell$, which partition each second in $I_h$ into 4 equidistantly spaced time instances.

We associate the Bonn-Beuel junction in Figure~\ref{fig:laneBB} with the network model introduced in Section \ref{sec:trafficflowmodels} and assume that \eqref{eq:lwr}, \eqref{eq:fundamentaldiagram} and \eqref{eq:riemannsolver} are a suitable model of the measured dynamics at the junction. In particular, we make the assumption that the vehicle densities \eqref{eq:density} at the junction interact with each other in a way that allows for a prediction of the traffic flow at the network node. In the following sections we address the question if this interaction can be described by the coupling models from Section~\ref{sec:couplingmodels}.
  
\subsection{Coupling Delay}\label{sec:delay}
In this section we are concerned with the timing of the coupling in the data. Due to our modeling approach and the choice of the volumes, the empirical vehicle densities \eqref{eq:density} do not necessarily all interact at the same time. Instead we assume that there are fixed time delays $\tau_\ell^k\in \R$ such that the coupling relation can be described by a Riemann solver, of the form
\begin{equation}\label{eq:datadelaycoupling}
\mathcal{RS}\left( \rho_d^1(t + \tau^1_\ell), \rho_d^2(t + \tau^2_\ell), \rho_d^3(t + \tau^3_\ell) \right) \approx \left(  f_d^1(t + \tau^1_\ell),  f_d^2(t + \tau^2_\ell),  f_d^3(t + \tau^3_\ell) \right). 
\end{equation}
Here we assume $t+\tau^k_\ell \in I_\ell$ for $k=1,2,3$. Since the average velocity, which contributes to the coupling delay, varies with the time of day (clf. Table~\ref{tab:datasets}) we allow for different delays in each data set.

\begin{figure*}[t]
  \centering
  \includegraphics{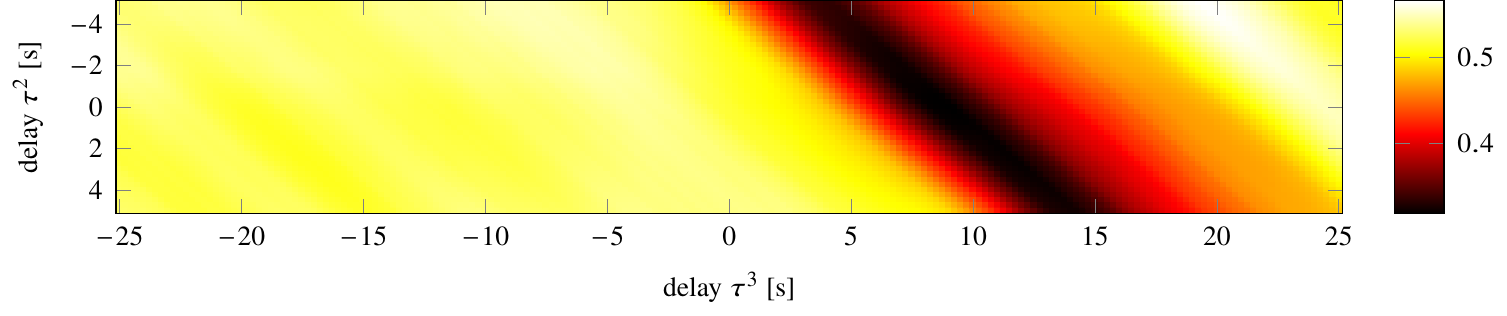}
  \vspace{5pt}
  \includegraphics{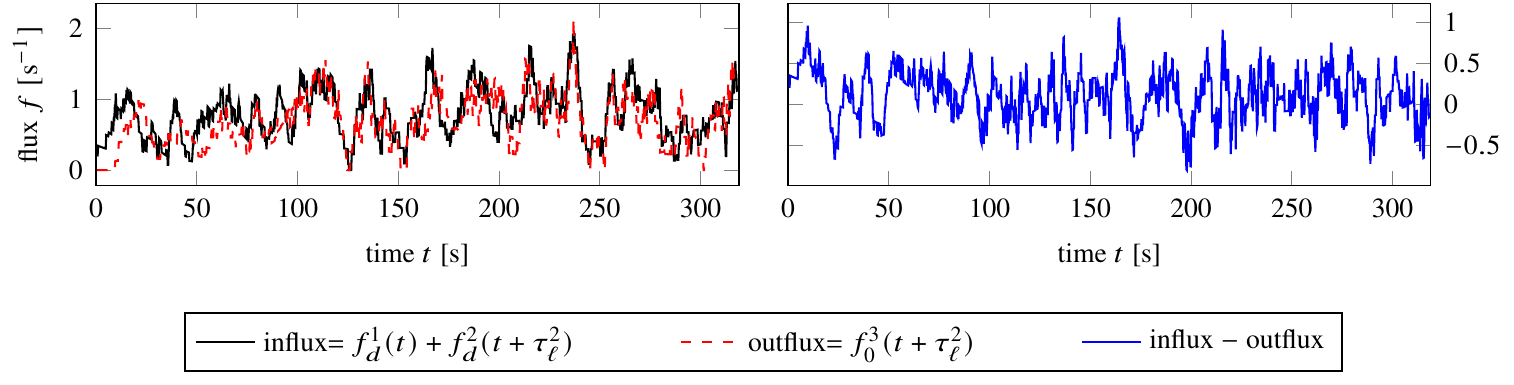}
  \caption{Estimation of the coupling delay in data set $\ell=13$. The optimization problem~\eqref{eq:delayoptimization} (top panel shows the functional to be minimized with respect to the delays) attains a minimum for $\tau^2=0$ and $\tau^3=9$. The accordingly shifted total influx approximately matches the shifted outflux qualitatively (bottom left) and the Kirchhoff condition is satisfied up to an absolute error of approximately $0.5$ point-wise in time (bottom right) on average.}\label{fig:delay}
\end{figure*}

To identify suitable delays, we set $\tau^1_\ell=0$ and assume that it holds $|\tau^2_\ell|\leq 5$ and $|\tau^3_\ell| \leq 25$ for all data sets $\ell=1,\dots,31$. We choose the delays such that the empirical traffic flow in \eqref{eq:datadelaycoupling} approximately satisfies the Kirchhoff condition \eqref{eq:kirchhoff}, i.e., as minimizers of the optimization problem
\begin{equation}\label{eq:delayoptimization}
  \min_{\tau^2_\ell, \tau^3_\ell} \| f_d^1(\cdot) + f_d^2(\cdot + \tau^2_\ell) - f_d^3(\cdot + \tau^2_\ell)\|_{I_\ell},
\end{equation}
for a norm taking into account the full time interval. We optimize over the discretized time intervals and use the root-mean-square norm
\begin{equation}\label{eq:rms}
  \| g \|_{I_\ell} = \sqrt{ \frac{1}{|I_\ell^h|} \sum_{t_i \in I^h_\ell} g(t_i)^2 } \quad \text{for} \quad g:I_\ell \rightarrow \R.
\end{equation}
In \eqref{eq:delayoptimization} the norm is taken over suitably shortened intervals such that $t+\tau_\ell^k \in I_\ell$ for $k=2,3$. For simplicity, we identify the minimizers by sampling, where we allow only for multiples of the grid size in $I_h$ (i.e., $0.25$) as candidates for the delays $\tau_\ell^2$ and $\tau_\ell^3$. 

In Figure~\ref{fig:delay} we show the optimization landscape of problem \eqref{eq:delayoptimization} in case of data set $\ell=13$. The unique global minimizer given by $\tau^2_\ell=0$ and $\tau^3_\ell=9$ determines the coupling delay in this case. The shifted empirical traffic flows approximately satisfy the Kirchhoff condition point-wise in time, see Figure~\ref{fig:delay} (bottom). The estimated coupling delays for all data sets are shown in Table~\ref{tab:datasets}. In general, similar delays were obtained at similar times of day. Positive $\tau^3$ were estimated for all data sets. The estimated $\tau^2$ range from $-5$ to $5$ and seem to be related to the difference in average speed between passing and entering vehicles.

In the rest of this work, we only consider shifted empirical quantities taking into account the coupling delays and, by abuse of notation, associate $\rho_d^k(t)$, $v_d^k(t)$ and $f_d^k(t)$ with $\rho_d^k(t+\tau_\ell^k)$, $v_d^k(t + \tau_\ell^k)$ and $f_d^k(t + \tau_\ell^k)$ assuming $t + \tau_\ell^k \in I_\ell$ for $k=1,2,3$. Moreover, we accordingly shorten the time intervals $I_\ell$ such that $t + \tau_\ell^k$ is included in the original time interval for all $t\in I_\ell$ and $k=1,2,3$.

\subsection{Fundamental Diagrams}\label{sec:fd}
We randomly grouped the 31 data sets into training data (consisting of 8 data sets), testing data (consisting also of 8 data sets) and application data (consisting of 15 data sets), for details see Appendix~\ref{sec:datasets}. The training data is used for parameter estimation for both the fundamental diagrams and the coupling models. While the testing data is used for first validations of the coupling models and to avoid over-fitting, the application data is employed to compare coupling models to traffic data on the full network. As the fundamental diagrams are employed in all considered coupling models they were estimated in the first step. 

\begin{figure}[t]
  \centering
  \includegraphics{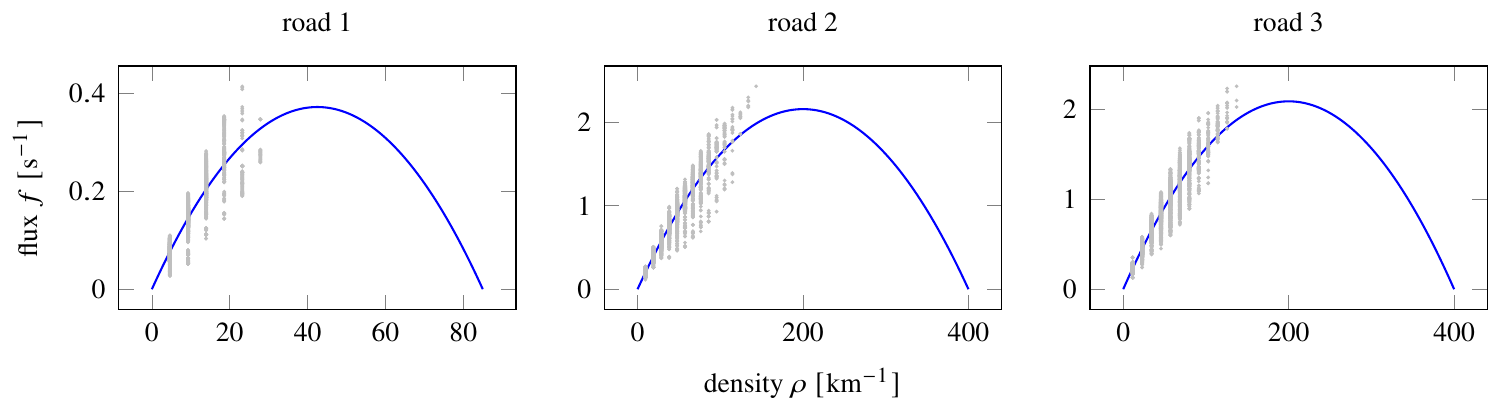}
  \caption{Fundamental diagrams of a 2-to-1 network corresponding to the Bonn-Beuel on-ramp (see Figure~\ref{fig:laneBB}). The density-flow form of the diagrams (blue lines) and the training data (gray dots) are shown. While the training data is congestion free when considering the individual roads congestion of the junction, i.e., $d^1+d^2>s^3$ occurs in the data.}\label{fig:FD}
\end{figure}

We fit the networked LWR model \eqref{eq:lwr} to the data by estimating the parameters of the fundamental diagram \eqref{eq:fundamentaldiagram} from the training data. We allow for different fundamental diagrams on the roads of the network but assume that they do not vary with respect to the time of day or the data set, respectively. Variability with the time of day in our approach is only reflected in the coupling delays. We use the empirical densities and velocities in the control volumes to estimate the fundamental diagram on the corresponding road of the network. We take a least-square approach and on each road $k=1,2,3$ select the parameters minimizing
\begin{equation}\label{eq:leastsquares}
  \sum_{t_i\in I^h_\text{train}} \left( v^k_\text{max}- \frac{v^k_\text{max}}{\rho^k_\text{max}} \rho^k_d(t_i) - v^k_d(t_i) \right)^2,
\end{equation}
where $I^{h,k}_\text{train}$ is the union of all time intervals considered in the training data excluding all $t_i$ with $\rho^k_d(t_i)=0$. As is often the case even for much larger data collections, see \cite{fan2014comparawrasclzhang}, large vehicle densities are not well represented in the data, which complicates the estimation of the stagnation densities of the fundamental diagrams. For this reason we constrain the parameters while minimizing \eqref{eq:leastsquares} using a typical stagnation density as upper bound, i.e.,
\begin{equation}\label{eq:stagmax} 
  \rho^k_\text{max} \leq \frac{\operatorname{card}\{ \text{lanes on road }k\}}{ 7.5 \text{ m}},
\end{equation}
where 7.5 m is an estimate of the minimal distance between the centers of vehicles taking the average vehicle length and safety distance into account. 
As result we obtain the maximal velocities $v^1_\text{max} = 62.94$ km/h, $v^2_\text{max} = 77.59$ km/h and $v^3_\text{max} = 75.28$ \kmh~and the stagnation densities $\rho^1_\text{max} = 84.99$ \pkm, $\rho^2_\text{max} = 400$ \pkm~and $\rho^3_\text{max} = 400$~\pkm. On roads 2 and 3 the maximal density due to~\eqref{eq:stagmax} is attained. We show the estimated fundamental diagrams together with the data in Figure~\ref{fig:FD}. While the training data does not include vehicle densities in the congestive regime with respect to the individual roads, congestion of the junctions, i.e. $d^1+d^2>s^3$ is observed in the training data (for approximately $0.17 \%$ of the data points).

\section{Model Fitting}\label{sec:fitting}
To calibrate the coupling models we estimate their parameters so that the model error with respect to the training data is minimized. The model error compares the empirical fluxes $f^k_d$ introduced in Section~\ref{sec:data} to the fluxes obtained by the coupling models applied to the empirical vehicle densities \eqref{eq:density} over the relevant time windows. Thus, the road dependent model errors read
\begin{equation}\label{eq:moderrorroad}
\mathcal{E}_k^M = \|f_d^k - \mathcal{RS}_k^M(\rho^1_d, \rho^2_d, \rho^3_d) \|_{I^h_\text{train}}^2
  \end{equation}
  for $k=1,2,3$, where $I^h_\text{train}$ denotes the union of all discretized time intervals for the training data and $\mathcal{RS}_k^M$ refers to the coupling flux for road $k$ corresponding to model~$M$. The root-mean-square norm \eqref{eq:rms} is used to compute an average over all time data. The average over the errors \eqref{eq:moderrorroad} for all coupling fluxes, i.e., influxes 1 and 2 and the outflux at the network node defines the total model error $\mathcal{E}^M$.
 
  Coupling models and parameter estimation are implemented in the Julia programming language \cite{bezanson2017julia}, we use the Flux library \cite{innes2018flux} to implement the machine learning models and the global optimization library~\cite{feldt2018black} for parameter estimation in the classical models. Our developed codes are publicly available from the repository~\cite{kolbe2022implemdata}.
  As we aim to estimate parameters over generally large and partly high-dimensional spaces and expect high nonlinearity of the coupling models we employ the differential evolution method~\cite{storn1997differ} for the classical models. This metaheuristics algorithm samples initial populations of parameter states that evolve by forming combinations between its members potentially replacing previous members if the combinations achieve a better fit. We use an adaptive variant of the original algorithm that initially guesses suitable control parameters and during the iteration optimizes them~\cite{qin2005self}. Although convergence of this algorithm has not been proved in our application it reliably identified minimizing parameter states in a feasible number of iterations.

  The classical coupling models C1--C4 are calibrated by minimizing the total model error $\mathcal{E}^M$ with respect to the training data over the corresponding parameter space. We allow the right of way parameter $\beta$ to take values within the full interval $[0,1]$ and constrain the empty road velocity in the second order models by $0.1\, v^k_\text{max} \leq w^k \leq v^k_\text{max}$ with $v^k_\text{max}$ taken from the estimation of the fundamental diagrams, see Section~\ref{sec:fd}. The computed optimal $\beta$ varies with the flow maximization approach; the coupling models using the second approach predicted $\beta\approx 5.23 \%$ in C3 and $\beta\approx 5.59 \%$ in C4 and the coupling models using the first approach predicted slightly smaller priorities ($\beta\approx 4.94 \%$ in C1 and $\beta\approx 3.11 \%$ in C2).
  
  \begin{table}\scriptsize
    \centering
    \caption{Maximal velocity on the roads due to the fundamental diagrams as computed in Section~\ref{sec:fd} (first column) and empty road velocities at the junction in the second order coupling models after parameter estimation.}\label{tab:vel}
    \begin{tabular}{ lcccc}
      \toprule
     & FD & C2 & C3 & C4  \\ \midrule 
      influx 1 & 62.94 km/h & 54.20 km/h & 60.55 km/h & 62.91 km/h \\
      influx 2 & 77.59 km/h & 71.28 km/h & 72.48 km/h & 74.07 km/h \\
      outflux & 75.28 km/h & 68.93 km/h & 73.18 km/h & 60.29 km/h \\ \bottomrule
    \end{tabular}
  \end{table}

  In Table~\ref{tab:vel} we compare the maximal velocities on the roads due to the fundamental diagrams in Section~\ref{sec:fd} to the empty road velocity at the junction in the second order coupling models C2--C4. While the velocities of the coupling models remain within a 10 \% interval around the velocity of the fundamental diagram in case of transit traffic on road 2, larger deviations on the entry lane (road 1) occur, most notable in model C2, which predicts an approximately 14 \% smaller velocity than the fundamental diagram. On the outgoing segment (road 3) model C4 assumes an exceptional decrease in velocity of approximately 20 \%  at the junction compared to the velocity on the road. 

  A plethora of training methods specialized for neural networks has been introduced in recent years. Instead of the metaheuristic algorithm used for the calibration of the classical models we use the Adam method \cite{kingma2014adam}, an extended version of stochastic gradient descent (SGD), to train the machine learning models introduced in Section~\ref{sec:mlmodels}. SGD approximates the descent direction by computing the gradient of the cost function based on a single randomly selected sample from the training data set. A series of successive SGD updates taking into account all training samples in this way is referred to as a training epoch. The Adam method makes use of running averages over the gradients and its second moments in the computation of the descent direction and to adapt the learning rate in each iteration. We specifically use the AMSGrad variant of the Adam method \cite{reddi2018converadambeyon}. To account for the constraint in the training of ML3 in \eqref{eq:consistencyconstraints} we introduce the penalty term
  \begin{equation}\label{eq:fdpenalty}
    \mathcal{P}_\text{C} = \left\| (f_0^1, f_0^2, f_0^3) - \mathcal{ML}( \rho_0^1, \rho_0^2, \rho_0^3)\right\|_{I^h_\text{train}}^2 
  \end{equation}
  and train the model using the cost function $\mathcal{E}^M + \frac12 \mathcal{P}_\text{C}$.

\begin{table}\scriptsize
  \caption{Training progress while testing the machine learning model capability to approximate the flow maximization model C1'. The mean and standard deviation over five training runs is shown.}\label{tab:loss}
  \centering
  \begin{tabular}{lrccc}
    \toprule
&   training epoch & ML1 & ML2 & ML3\\ \midrule
\multirow{5}{*}{\rotatebox[origin=c]{90}{training loss}} & 0 & \(2.897 \times 10^{-3} \pm 1.141 \times 10^{-3}\) & \(3.618 \times 10^{-3} \pm 1.795 \times 10^{-3}\)  & \(5.152 \times 10^{-3} \pm 1.328 \times 10^{-3}\)\\[0pt]
 & 1 & \(4.154 \times 10^{-4} \pm 3.250 \times 10^{-5}\) & \(5.285 \times 10^{-4} \pm 2.068 \times 10^{-5}\) & \(2.037 \times 10^{-3} \pm 2.934 \times 10^{-5}\)\\[0pt]
& 10 & \(2.365 \times 10^{-4} \pm 7.822 \times 10^{-6}\) & \(2.362 \times 10^{-4} \pm 2.619 \times 10^{-5}\) & \(9.552 \times 10^{-4} \pm 2.239 \times 10^{-5}\)\\[0pt]
& 100  & \(1.935 \times 10^{-4} \pm 2.316 \times 10^{-7}\) & \(1.056 \times 10^{-5} \pm 9.123 \times 10^{-7}\) & \(3.609 \times 10^{-4} \pm 1.416 \times 10^{-5}\)\\[0pt]
& 500 & \(1.931 \times 10^{-4} \pm 2.909 \times 10^{-7}\) & \(6.516 \times 10^{-6} \pm 5.393 \times 10^{-7}\) & \(1.127 \times 10^{-4} \pm 1.736 \times 10^{-5}\)\\[5pt]
\multirow{5}{*}{\rotatebox[origin=c]{90}{test loss}}& 0  & \(2.811 \times 10^{-3} \pm 1.134 \times 10^{-3}\) & \(3.508 \times 10^{-3} \pm 1.804 \times 10^{-3}\) & \(5.309 \times 10^{-3} \pm 1.268 \times 10^{-3}\)\\[0pt]
& 1 & \(3.712 \times 10^{-4} \pm 2.730 \times 10^{-5}\) & \(4.573 \times 10^{-4} \pm 2.536 \times 10^{-5}\) & \(2.416 \times 10^{-3} \pm 4.301 \times 10^{-5}\)\\[0pt]
& 10 & \(2.136 \times 10^{-4} \pm 7.838 \times 10^{-6}\) & \(2.149 \times 10^{-4} \pm 3.025 \times 10^{-5}\) & \(8.923 \times 10^{-4} \pm 2.154 \times 10^{-5}\)\\[0pt]
& 100 & \(1.711 \times 10^{-4} \pm 8.832 \times 10^{-7}\) & \(1.079 \times 10^{-5} \pm 9.557 \times 10^{-7}\) & \(3.625 \times 10^{-4} \pm 1.907 \times 10^{-5}\)\\[0pt]
& 500 & \(1.705 \times 10^{-4} \pm 5.025 \times 10^{-7}\) & \(7.570 \times 10^{-6} \pm 6.491 \times 10^{-7}\) & \(1.578 \times 10^{-4} \pm 2.031 \times 10^{-5}\)\\ \bottomrule
\end{tabular}
\end{table}

To test the capabilities of the machine learning model we use the model architecture of ML1--ML3 to approximate/learn the flow maximization model in a standardized setting. Let C1' refer to the flow maximization model as defined in Section~\ref{sec:classic} assuming that the velocities on the roads are given by the Greenshield model \eqref{eq:fundamentaldiagram} with parameters chosen $v_\text{max}^k=\rho_\text{max}^k=1$ for $k\in\{1,2,3\}$ and priority set to $\beta=0.5$. Discretizing the set of suitable input traces $[0,1]^3$ by 20 equidistantly spaced points in each of the three dimensions and computing the coupling fluxes using C1' defines the training data set. The test data set is similarly defined using a finer discretization by 80 points per dimension. We train the three machine learning models over 500 epochs using this training data and present the mean and standard deviation of both the training and the test loss over five training runs in Table~\ref{tab:loss}. In case of ML3 the penalty term is neglected in the table. In all machine learning models the training loss decreases over the training. While in ML1 the decrease mostly stagnates after 100 epochs further significantly improvements in ML2 and ML3 with respect to the training loss are still visible in later epochs. The model ML2 is trained faster than ML1 and ML3 and achieves a training loss of approximately  $6 \times 10^{-6}$ after 500 epochs. In all models the test loss continuously follows the decrease of the training loss. Experiments with different architectures have revealed a significant benefit of ML3 over two and three layer approaches with respect to the training performance, whereas improvements have become less significant when even higher number of layers have been considered.

We analogously train our machine learning models to the data from Section~\ref{sec:data}. Although it minimizes the costs only locally the Adam method significantly reduces the costs in its first training era. Testing against the test data to avoid over-fitting training is conducted over 100 training eras in total.

\begin{table}\scriptsize
  \caption{Road dependent and total mean square model errors \eqref{eq:moderrorroad} with respect to the training and the test data for the introduced coupling models.}\label{tab:moderrors}
  \centering
    \begin{tabular}{lcccccccc}
      \toprule
      &  & C1 & C2 & C3 & C4 & ML1 & ML2 & ML3\\ \midrule
      \multirow{4}{*}{\rotatebox[origin=c]{90}{training data}} &  influx 1& \(4.816 \times 10^{-4}\) & \(7.612 \times 10^{-4}\) & \(9.122 \times 10^{-3}\) & \(8.693 \times 10^{-3}\) & \(4.470 \times 10^{-3}\) & \(3.456 \times 10^{-3}\) & \(3.642 \times 10^{-3}\)\\[0pt]
 &  influx 2 & \(1.091 \times 10^{-2}\) & \(1.687 \times 10^{-2}\) & \(1.880 \times 10^{-1}\) & \(1.871 \times 10^{-1}\) & \(2.498 \times 10^{-2}\) & \(2.595 \times 10^{-2}\) & \(2.418 \times 10^{-2}\)\\[0pt]
 &  outflux & \(1.843 \times 10^{-1}\) & \(1.667 \times 10^{-1}\) & \(2.773 \times 10^{-1}\) & \(2.782 \times 10^{-1}\) & \(1.204 \times 10^{-1}\) & \(1.167 \times 10^{-1}\) & \(1.182 \times 10^{-1}\)\\[0pt]
 &  total & \(6.523 \times 10^{-2}\) & \(6.143 \times 10^{-2}\) & \(1.581 \times 10^{-1}\) & \(1.580 \times 10^{-1}\) & \(4.996 \times 10^{-2}\) & \(4.871 \times 10^{-2}\) & \(4.866 \times 10^{-2}\)\\[5pt]
 \multirow{4}{*}{\rotatebox[origin=c]{90}{test data}} & influx 1 & \(2.309 \times 10^{-4}\) & \(6.723 \times 10^{-4}\) & \(8.040 \times 10^{-3}\) & \(7.668 \times 10^{-3}\) & \(3.523 \times 10^{-3}\) & \(2.954 \times 10^{-3}\) & \(3.155 \times 10^{-3}\)\\[0pt]
 & influx 2 & \(6.226 \times 10^{-3}\) & \(1.505 \times 10^{-2}\) & \(2.830 \times 10^{-1}\) & \(2.817 \times 10^{-1}\) & \(1.986 \times 10^{-2}\) & \(2.132 \times 10^{-2}\) & \(2.030 \times 10^{-2}\)\\[0pt]
 & outflux & \(7.921 \times 10^{-2}\) & \(6.651 \times 10^{-2}\) & \(2.437 \times 10^{-1}\) & \(2.440 \times 10^{-1}\) & \(3.693 \times 10^{-2}\) & \(3.256 \times 10^{-2}\) & \(3.385 \times 10^{-2}\)\\[0pt]
 & total & \(2.856 \times 10^{-2}\) & \(2.741 \times 10^{-2}\) & \(1.782 \times 10^{-1}\) & \(1.778 \times 10^{-1}\) & \(2.010 \times 10^{-2}\) & \(1.894 \times 10^{-2}\) & \(1.910 \times 10^{-2}\)\\ \bottomrule
  \end{tabular}
  \end{table}

In Table~\ref{tab:moderrors} we present the road dependent and total model errors of all the coupling models for best-fit parameters considering both the training and the test data. Among the classical models C2 achieve the best total fit to the data with its first order variant C1 achieving only slightly larger errors. The second order models C3 and C4 achieve the largest total errors within all tested models. With respect to the total model error the machine learning models outperform the classical ones. All three yield similar total model errors with ML3 yielding the smallest with respect to the training data and ML2 the smallest with respect to the test data. Overall the results are consistent both in the training and in the test data. Surprisingly, most models (all except for C3 and C4) yield an even better fit to the test data than to the training data. This might be due to the fact that the training data covers more situations of junction congestion compared to the test data, see Appendix~\ref{sec:datasets}.
Regarding the roads at the junction all models achieve the best fit in the coupling flux with respect to the entry lane with the classical models performing significantly better in the training data when compared to the machine learning models.

\begin{figure}[t]
  \centering
  \includegraphics{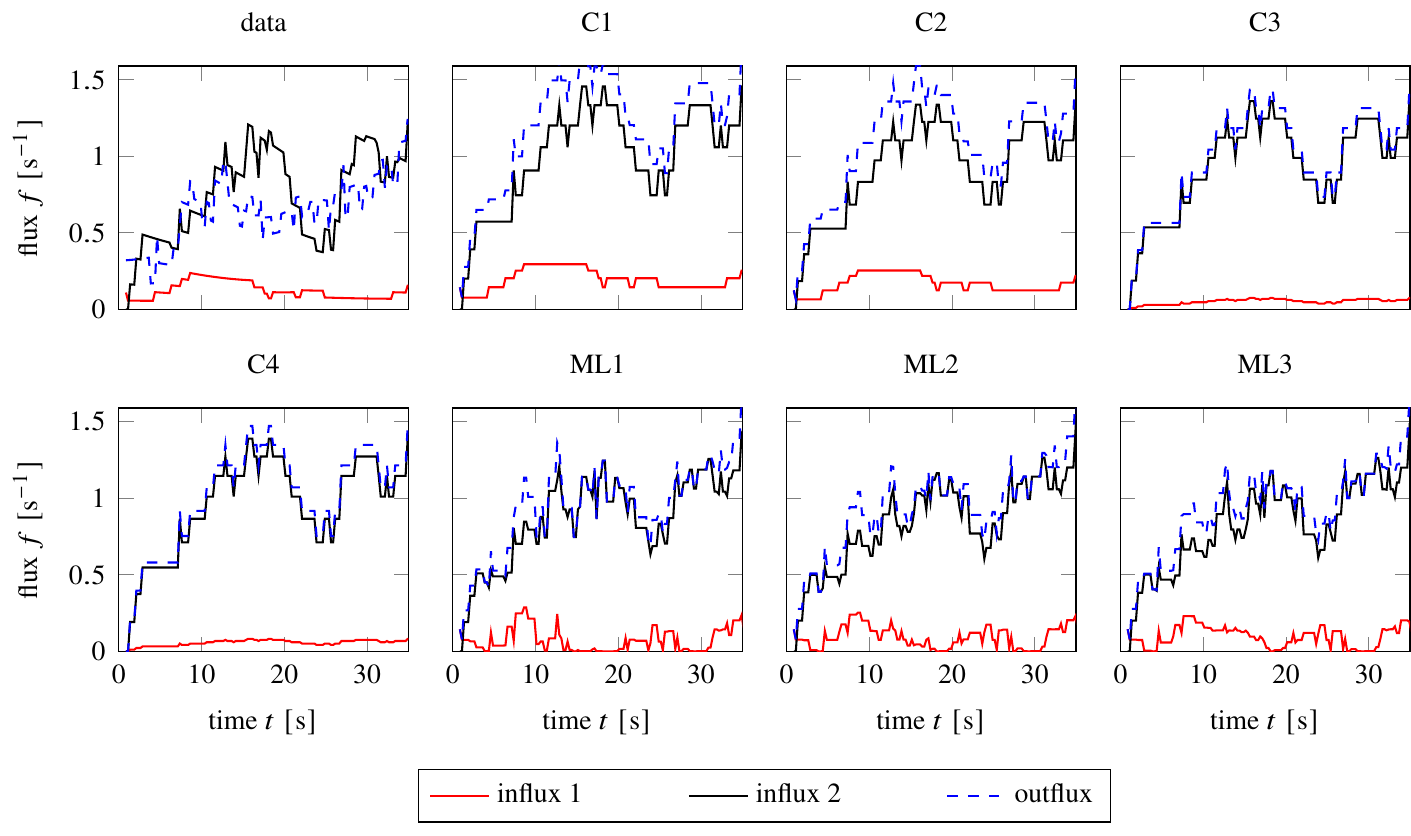}
  \caption{Empirical fluxes (top left) compared to coupling flux predictions by the coupling models. A time window with busy traffic from data set 10 within the training data set (see Table \ref{tab:datasets}) is considered.}\label{fig:fluxes}
\end{figure}

Figure~\ref{fig:fluxes} provides a qualitative insight into the flux predictions of the coupling models. We consider approximately half a minute of busy traffic on the Bonn-Beuel on-ramp captured in the training data. The time evolution of the empirical fluxes is shown in the upper left panel. It is evident that during the observed time window the Kirchhoff condition \eqref{eq:kirchhoff} is partly violated as the influx from road 2 (and therefore clearly the total influx) exceeds the outflux. This is due to the approach we use to compute macroscopic data not being able to exactly capture the coupling fluxes at the junction. The coupling models correct this as they all predict the outflux to exceed the influx from road 2 at almost all times. The models C1 and C2 generally predict the largest differences between the concerned fluxes. Both predict an increased influx from road 2 compared to the data with qualitatively similar evolution in time. They slightly overestimate the influx from the on-ramp in comparison to the data. Both models exhibit almost indistinguishable fluxes with only model C1 yielding sometimes larger outflux. Also, the results of models C3 and C4 are very similar to each other. In comparison to models C1 and C2 they predict smaller outflux and influx from the on-ramp. While the fluxes by the classical models, similar to the data, develop plateaus, i.e., constant fluxes on a shorter time-scale, the fluxes by the machine learning models appear more time varying and almost noisy. Unlike the classical models the machine learning ones predict influx from the on-ramp around zero at certain times. All three exhibit again almost similar dynamics but minor differences between the results of the models; an example occurs  around time $t=15$, where the single layer models predict almost zero influx from the on-ramp which is not the case in the four layer models ML2 and ML3. 

\section{Model Validation}\label{sec:validation}
In this section we validate the developed coupling fluxes by investigating their behavior in the setting of the network model~\eqref{eq:lwr}. To this end we associate the freeway junction on Figure~\ref{fig:laneBB} with the network in Figure~\ref{fig:2to1}. We assume left boundaries of road 1 and 2 with x-position at the left interface of volume $V_2$ and a right boundary of road 3 with x-position at the right edge of volume~$V_3$. Considering these boundaries we estimate the length of the junction to be $2s=270.28$~m and for simplicity assume that all roads of the network model are of length $s$. On the network we impose the fundamental diagram~\eqref{eq:fundamentaldiagram} with road dependent parameters as estimated in Section \ref{sec:fd}.

In order to test the coupling models and compare the vehicle densities of the network model a numerical scheme is required. We consider an adapted variant of a recently introduced central scheme for networks \cite{herty2022centr}. The roads and their parametrisations $[-s,0]$ (road 1 and road 2) and $[0, s]$ (road 3) are discretized in $m$ equidistant cells of length $\Delta x$. Let $\rho_j^{k,n}$ denote an approximate average of $\rho^k$ over the cell $I_j=[(j-1/2)\Delta x, (j+1/2)\Delta x]$ at time instance $t_n$ and consider a time increment $\Delta t$. The scheme can then be written in conservative form as
\begin{equation}\label{eq:conservativeformnet}
  \rho_j^{k,n+1} = \rho_j^{k,n} - \frac{\Delta t }{\Delta x} \left( F_{j+1/2}^{k,n} - F_{j-1/2}^{k,n}\right).
\end{equation}
 In our numerical experiments we choose the time increment according to the CFL condition $\Delta t = \text{CFL}\, \frac{\Delta x}{\lambda}$. The cell indicating index in \eqref{eq:conservativeformnet} can be taken $j= -m,\dots, -1$ for $k=1,2$ or $j= 0, 1, \dots, m$ for $k=3$. The numerical fluxes are given by
\begin{equation}\label{eq:netfluxes}
  F_{j-1/2}^{k,n} =
  \begin{cases}
    \frac 1 2 \, (f^k(\rho_{j}^{k,n}) + f^k(\rho_{j-1}^{k,n}))  - \frac {\lambda} 2 (\rho_j^{k,n}-\rho_{j-1}^{k,n}) & \text{if }j \neq 0,\\[5pt]
    \mathcal{RS}^M_k(\rho^{1,n}_{-1}, \rho^{2,n}_{-1}, \rho^{3,n}_0) &\text{if }j=0.
  \end{cases}
\end{equation}
The coupling node of the network is located at the interface between the cells $C_{-1}$ and $C_0$ and therefore the coupling models give rise to the fluxes $F_{-1/2}^{k,n}$. The relaxation speed $\lambda>0$ is a parameter that stems from a relaxation approach used in the scheme derivation, see \cite{jin1995relaxschemsystem}. It is chosen such that the subcharacteristic condition $-\lambda\leq f_k^\prime(\rho^k) \leq \lambda$ holds on all roads. To also account for the coupling fluxes we adaptively update $\lambda$ (and consequently also the time increment) in each time step such that
\[
  \lambda> \frac{2}{\Delta x}\min \left\{  \left| F_{-1/2}^{1,n}- f^1(\rho_{-1}^{k,n}) \right|,~\left| F_{-1/2}^{2,n}- f^2(\rho_{-1}^{k,n}) \right|,~\left|  f^k(\rho_{0}^{3,n}) - F_{-1/2}^{3,n}\right| \right\}.
\]

\subsection{Boundary Fluxes}\label{sec:bfluxes}
In this section we derive ingoing and outgoing fluxes at the network boundaries from the application data and use them to validate the coupling fluxes. Employing the freeway junction boundaries as defined above, we can for any vehicle trajectory $\mathbf x_j$ compute the time at which the corresponding vehicle enters the network on either road 1 or road 2 and the time at which it leaves the network again from road 3. If a vehicle is already on the junction at start time of the recording or has not left the network at end time of the recording the times of entering and leaving can be computed by extrapolating its polynomial trajectory. Taking into account a full dataset $\mathcal{D}_\ell$ the entering and leaving times are used to compute the histograms
\begin{equation}
  V^k(t) = \operatorname{card}\{ \text{vehicles passing the boundary on road }k \text{ between time } \lfloor t \rfloor \text{ and }  \lfloor t \rfloor + 1\} 
\end{equation}
for $t\in I_\ell$ and $k=1,2,3$. These histograms are a representation of the boundary fluxes in the data. As the discontinuity of $V$ with respect to time might lead to instabilities when applied in the time- and space continuous network model \eqref{eq:lwr} we consider a regularization using a (weighted) kernel density estimator with Gaussian kernel and define
\begin{equation}
  \hat{V}^k(t) = \frac{1}{h} \sum_{t_i \in \mathcal{T}^k_\ell} \Phi \left( \frac{t- t_i}{h}\right),
\end{equation}
where the set $\mathcal{T}^k_\ell$ includes all entering/leaving times on road $k$ within the data set $\mathcal{D}_\ell$, the bandwith is chosen $h=0.75$ and $\Phi$ denotes the standard normal density function.

Fixing an application data set $\mathcal{D}_\ell$ we use the incoming boundary fluxes from the data as boundary information for the network. In the scheme given by~\eqref{eq:conservativeformnet} and \eqref{eq:netfluxes} we therefore impose the left boundary fluxes
\begin{equation}\label{eq:leftfluxes}
F_{-m-1/2}^{k, n} = \frac 1 2 (f^k(\rho_{-m}^{k,n}) + \hat{V}^k(t^n)) \quad \text{for }k=1,2
\end{equation}
in correspondence to the interior fluxes in \eqref{eq:netfluxes}. On the right boundary we instead impose homogeneous Neumann data, i.e. $\ddx \rho^3(s) =0$, and keep track of fluxes going out of the network, which are given by $V^3_{M}(t^n)= f^3(\rho^{3,n}_m)$. Initial data is taken constantly zero. In the computations we use $m=200$ cells per road, the Courant number $\text{CFL}=0.24$ and minimal relaxation speed $\lambda=10$.

\begin{figure}[t]
  \centering
  \includegraphics{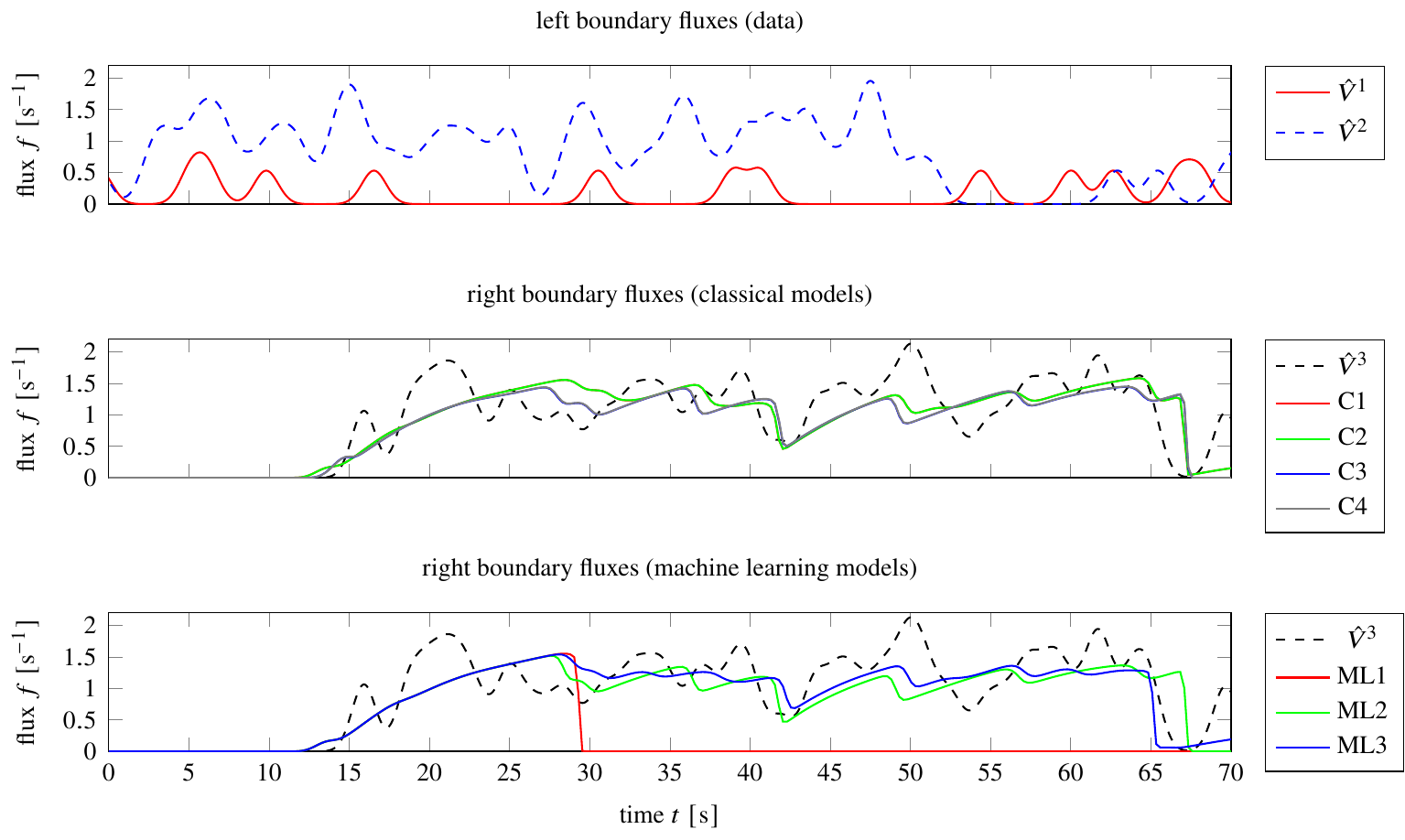}
  \caption{Left and right boundary fluxes over time from dataset 18 and model predictions by the network for classical (center) and machine learning (bottom) coupling models. The predictions at the right boundary are simulated using the left boundary fluxes from the data (top) as boundary conditions. The models C3 and C4 as well as ML1 and ML2 predict indistinguishable results, respectively.}\label{fig:boundaryfluxes}
\end{figure}

We present a time excerpt of the boundary fluxes from the data set 18 in Figure~\ref{fig:boundaryfluxes} together with the model predictions $V^3_{M}(t^n)$. While the models cannot reproduce the exact structure of the right boundary fluxes from the data, the classical models predict the boundary outflow reasonably well. In particular, we observe a good indication at which times boundary outflow occurs and also the position of some extrema of the trajectory is fit well. The models C1 and C2 predict a slightly larger magnitude of the trajectory and are therefore closer to the data than the models C3 and C4. For the given boundary inflow the one layer model ML1 apparently block the junction once a certain density threshold is reached and therefore the predicted right boundary flux first rises but then falls to zero after 30 seconds and afterwards stays constant; therefore this models do not fit the data. The four layer models ML2 and ML3 perform significantly better, they show a similar qualitative behavior as the models C3 and C4 with ML3 predicting generally larger fluxes and a sharper drop around $t=65$ due to the gap of passing cars. In Table~\ref{tab:boundaryerrors} we present the relative errors
\begin{equation}\label{eq:berror}
  \frac{\|\hat{V}^3 - V^3_{M} \|_{L^2(I_\ell)}}{\|\hat{V}^3\|_{L^2(I_\ell)}}
\end{equation}
with respect to all coupling models and data sets from the application data. The computed errors and their averages over all considered data sets generally confirm our findings from the discussion of Figure~\ref{fig:boundaryfluxes} and show that ML2 performs similarly well as the classical models C1 and C2 while the relative error of ML3 is in a similar range to the ones of C3 and C4.

\begin{table}\fontsize{8}{8}\selectfont
  \centering
    \caption{Relative $L^2$ errors \eqref{eq:berror} of the predicted boundary fluxes with respect to the data. Errors for all coupling models and all data sets from the application data are shown.}\label{tab:boundaryerrors}
  \begin{tabular}{rccccccc}
    \toprule
    \multicolumn{1}{c}{data set} & C1 & C2 & C3 & C4 &ML1 &ML2 & ML3 \\ \midrule
    2 & \(4.715 \times 10^{-1}\) & \(4.708 \times 10^{-1}\) & \(4.937 \times 10^{-1}\) & \(4.935 \times 10^{-1}\) & \(9.281 \times 10^{-1}\) & \(4.866 \times 10^{-1}\) & \(5.660 \times 10^{-1}\)\\[0pt]
3 & \(4.404 \times 10^{-1}\) & \(4.385 \times 10^{-1}\) & \(4.463 \times 10^{-1}\) & \(4.458 \times 10^{-1}\) & \(9.818 \times 10^{-1}\) & \(4.558 \times 10^{-1}\) & \(6.118 \times 10^{-1}\)\\[0pt]
5 & \(3.689 \times 10^{-1}\) & \(3.691 \times 10^{-1}\) & \(3.993 \times 10^{-1}\) & \(3.984 \times 10^{-1}\) & \(9.952 \times 10^{-1}\) & \(3.979 \times 10^{-1}\) & \(3.851 \times 10^{-1}\)\\[0pt]
7 & \(3.557 \times 10^{-1}\) & \(3.559 \times 10^{-1}\) & \(4.579 \times 10^{-1}\) & \(4.632 \times 10^{-1}\) & \(4.187 \times 10^{-1}\) & \(3.584 \times 10^{-1}\) & \(3.494 \times 10^{-1}\)\\[0pt]
9 & \(4.355 \times 10^{-1}\) & \(4.369 \times 10^{-1}\) & \(4.657 \times 10^{-1}\) & \(4.651 \times 10^{-1}\) & \(9.544 \times 10^{-1}\) & \(4.696 \times 10^{-1}\) & \(6.854 \times 10^{-1}\)\\[0pt]
12 & \(5.226 \times 10^{-1}\) & \(5.223 \times 10^{-1}\) & \(5.344 \times 10^{-1}\) & \(5.375 \times 10^{-1}\) & \(1.002\) & \(5.331 \times 10^{-1}\) & \(6.679 \times 10^{-1}\)\\[0pt]
13 & \(4.929 \times 10^{-1}\) & \(4.926 \times 10^{-1}\) & \(5.256 \times 10^{-1}\) & \(5.254 \times 10^{-1}\) & \(1.001\) & \(5.294 \times 10^{-1}\) & \(5.018 \times 10^{-1}\)\\[0pt]
17 & \(4.841 \times 10^{-1}\) & \(4.842 \times 10^{-1}\) & \(5.727 \times 10^{-1}\) & \(5.814 \times 10^{-1}\) & \(8.827 \times 10^{-1}\) & \(5.057 \times 10^{-1}\) & \(4.949 \times 10^{-1}\)\\[0pt]
18 & \(4.927 \times 10^{-1}\) & \(4.926 \times 10^{-1}\) & \(5.102 \times 10^{-1}\) & \(5.105 \times 10^{-1}\) & \(9.634 \times 10^{-1}\) & \(5.197 \times 10^{-1}\) & \(4.982 \times 10^{-1}\)\\[0pt]
23 & \(4.981 \times 10^{-1}\) & \(4.982 \times 10^{-1}\) & \(6.220 \times 10^{-1}\) & \(6.940 \times 10^{-1}\) & \(9.473 \times 10^{-1}\) & \(4.872 \times 10^{-1}\) & \(7.586 \times 10^{-1}\)\\[0pt]
25 & \(4.580 \times 10^{-1}\) & \(4.584 \times 10^{-1}\) & \(9.319 \times 10^{-1}\) & \(1.008\) & \(9.369 \times 10^{-1}\) & \(4.443 \times 10^{-1}\) & \(4.584 \times 10^{-1}\)\\[0pt]
27 & \(5.378 \times 10^{-1}\) & \(5.385 \times 10^{-1}\) & \(5.699 \times 10^{-1}\) & \(5.705 \times 10^{-1}\) & \(9.848 \times 10^{-1}\) & \(5.582 \times 10^{-1}\) & \(7.228 \times 10^{-1}\)\\[0pt]
28 & \(4.872 \times 10^{-1}\) & \(4.874 \times 10^{-1}\) & \(5.374 \times 10^{-1}\) & \(5.383 \times 10^{-1}\) & \(9.785 \times 10^{-1}\) & \(4.970 \times 10^{-1}\) & \(7.020 \times 10^{-1}\)\\[0pt]
29 & \(4.710 \times 10^{-1}\) & \(4.717 \times 10^{-1}\) & \(5.417 \times 10^{-1}\) & \(5.503 \times 10^{-1}\) & \(9.455 \times 10^{-1}\) & \(4.837 \times 10^{-1}\) & \(5.834 \times 10^{-1}\)\\[0pt]
31 & \(5.146 \times 10^{-1}\) & \(5.146 \times 10^{-1}\) & \(7.022 \times 10^{-1}\) & \(7.541 \times 10^{-1}\) & \(9.452 \times 10^{-1}\) & \(5.250 \times 10^{-1}\) & \(5.123 \times 10^{-1}\)\\ \midrule
average & \(4.687 \times 10^{-1}\) & \(4.688 \times 10^{-1}\) & \(5.541 \times 10^{-1}\) & \(5.691 \times 10^{-1}\) & \(9.244 \times 10^{-1}\) & \(4.834 \times 10^{-1}\) & \(5.665 \times 10^{-1}\)\\ \bottomrule
  \end{tabular}
  \end{table}

  \subsection{Model Prediction}
  In a last test we analyze the capability of the coupling models to predict dynamics not observed in the data. In particular, we consider congestion of the freeway. Again we employ the network model \eqref{eq:lwr} and the numerical scheme given by \eqref{eq:conservativeformnet} and \eqref{eq:netfluxes}. We employ the same scheme parameters as done in Section~\ref{sec:bfluxes} but this time assume homogeneous Neumann boundary conditions at all boundaries. We impose a Riemann problem with road-wise constant initial data
  \begin{equation}\label{eq:ics}
    \rho^{1,0} \equiv 0.7 \, \rho^1_\text{max}, \qquad \rho^{2,0} \equiv 0.5 \, \rho^2_\text{max}, \qquad \rho^{3,0} \equiv 0.8 \, \rho^3_\text{max},
  \end{equation}
  i.e., we assume congestion behind the on-ramp as the vehicle density reaches 80\% of the road capacity and busy traffic on the motorway before the junction (50\% occupied) as well as a congested entry lane (70\% occupied).
  
\begin{figure}[t]
  \centering
  \includegraphics{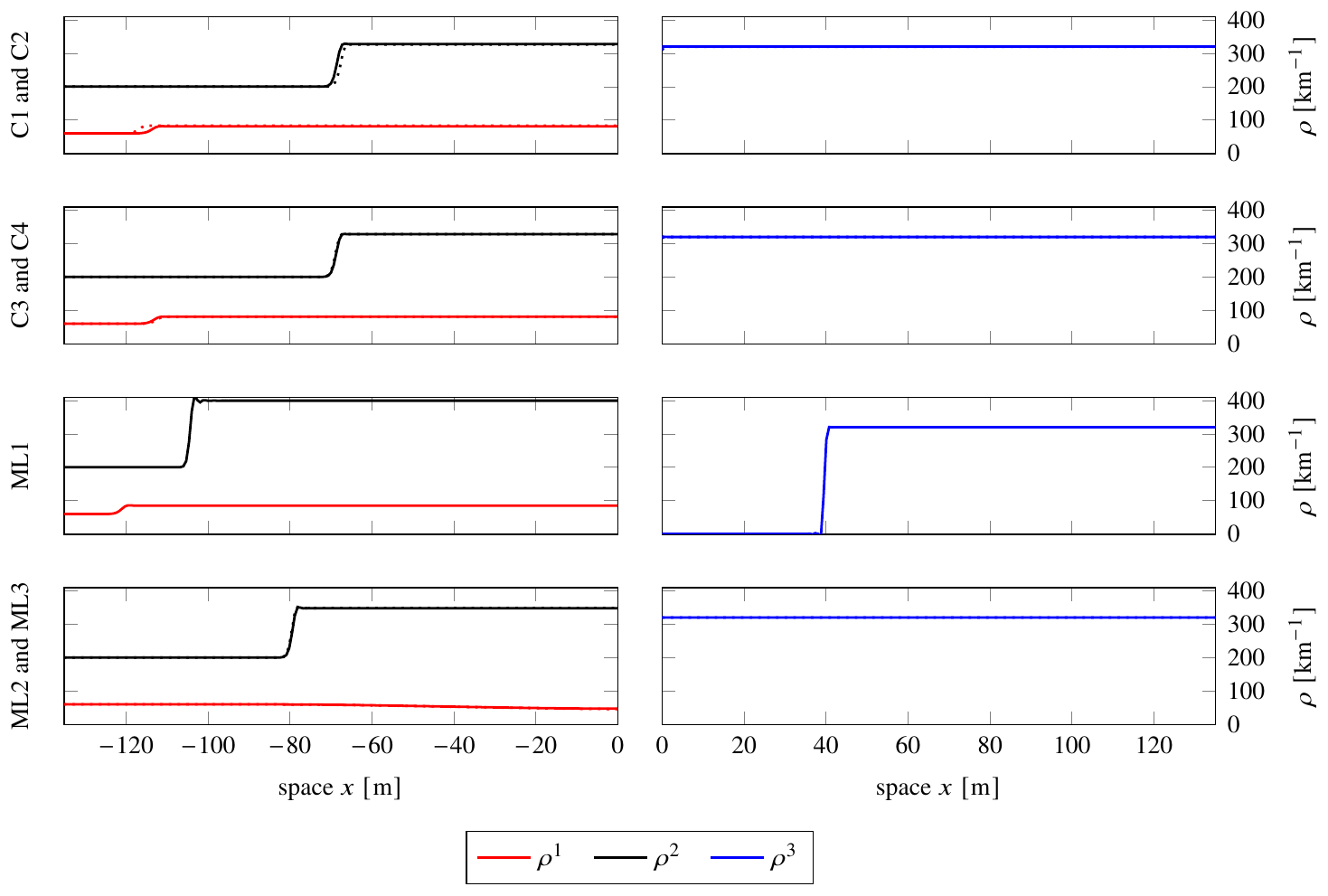}
  \caption{Numerical vehicle densities on the incoming roads 1 and 2 (left) and on the outgoing road 3 (right) of the traffic network as predicted 10 seconds after the congested state \eqref{eq:ics} by the coupling models. Dotted lines represent the respective second model of the panel. The solutions are computed by the scheme~\eqref{eq:conservativeformnet}.}\label{fig:prediction}
\end{figure}

Figure~\ref{fig:prediction} shows the model prediction addressing the evolution of this state after 10 seconds. Models C1 and C2 predict congestion on roads 1 and 2 as backwards traveling traffic waves on both roads occur, where the one on the entry-lane moves faster. Model C2 suggests a slightly slower propagation of the wave on the freeway compared to C1. In models C3 and C4 a slight reduction in wave speed over C1 is also observed on the entry-lane. As observed in Section~\ref{sec:bfluxes} the single layer model ML1 blocks the junction and no traffic is admitted to road 3. Congestion waves occur on road 1 and road 2 with magnitudes close to the stagnation densities and with the wave on road~2 being significantly faster than in the classical models. The models ML3 and ML4 predict a different behavior from the other models: a backward propagating congestion wave occurs only on road~2 while the vehicle density slightly decreases on the entry-lane. The congestion wave on road~2 is predicted faster than in the classical models. These dynamics are comparable to a situation, where the reduction of traffic congestion for entering cars is assigned a higher priority than the flow on the main lanes of the freeway.

\section{Conclusion}
We have presented a new modeling framework for junctions based on vehicle trajectory data using macroscopic models of traffic flow. Focusing, in particular, on the LWR model we have discussed the modeling of traffic networks and introduced a unified description of coupling models in terms of Riemann solvers. To preserve the mass of vehicles at the junction the models need to satisfy the Kirchhoff- as well as the demand and supply conditions. As these conditions are not sufficient to obtain well-defined coupling conditions we considered various coupling models that determine the coupling fluxes. We have introduced a new class of machine learning models that combine problem specific layers with ANNs. Employing a new parameterized representation of the admissible set the ANNs are trained to select suitable coupling fluxes.  These new models have been compared to established coupling rules that mainly rely on flow maximization. 
We have tested our framework using vehicle trajectory data from an on-ramp junction, and thereby compared the considered coupling models. Macroscopic data in terms of vehicle densities and velocities has been derived by tracking vehicles in pre-defined control volumes right before and after vehicles from the different roads interact with each other. To both account for the synchronization between the interaction and variability with respect to the time of day we propose a normalization based on time delays by fitting the data to the Kichhoff condition. The normalized macroscopic data is further used to estimate the fundamental diagram, relating vehicle densities and velocities, on each road. Employing training data covering traffic in the environment of the on-ramp of approximately 37 minutes the model parameters have been estimated. For this purpose an evolutionary algorithm and a variant of stochastic gradient descent in case of the neural network has been used.
The machine learning models have achieved a better fit to the data than the classical ones in both, the training data and an equally large test data set. The best-fit models have been tested in the networked PDE setting. Using additional data over the on-ramp we have compared throughput of the junction predicted by the coupling models combined with a finite volume scheme and the measurements. In a second test the models have been used to predict traffic in case of congestion, which has not been well represented in the data. The classical models based on flow maximization achieved a good fit to the data and reasonable predictions in these tests. While single layer machine learning models have still yield unrealistic network solution (complete blockage of the junction) an increase to four layers has resulted in a boundary flux prediction similarly accurate as the one obtained by the best classical models in the first test and realistic results (strong prioritization of the traffic from the on-ramp) significantly differing from the prediction of the classical models in the second one.
The results indicate that coupling models based on machine learning constitute a promising alternative to classical coupling rules as they are also able to predict realistic phenomena of traffic that are not present in the classical models. In addition, these models might further improve if more training data is available.

\subsection*{Acknowledgments}
The authors thank the Deutsche Forschungsgemeinschaft (DFG, German Research Foundation) for the financial support through HE5386/18-1,19-2,22-1,23-1. Furthermore, the authors would like to thank the Institute of Highway Engineering at RWTH Aachen University, in particular, M. Oeser, E. Kallo and M. Berghaus for kindly providing access to the trajectory data. For the evaluation of model C3 and C4, we thank S: G\"ottlich and J. Weissen (University of Mannheim) for providing the corresponding source code used in the paper \cite{göttlich2021seconordertraff}.

\appendix

\section{Data sets}\label{sec:datasets}
Table~\ref{tab:datasets} provides details on the considered data sets, such as time of day, number and average speed of passing and entering cars and the estimated coupling delays in Section~\ref{sec:delay}. Training data consists of the 8 data sets 1, 4, 10, 11, 16, 19, 20 and 30, test data of the 8 data sets 6, 8, 14, 15, 21, 22, 24 and 26 and application data of the 15 data sets 2, 3, 5, 7, 9, 12, 13, 17, 18, 23, 25, 27, 28, 29 and 31. Out of the two data sets capturing slowed down traffic (data sets 10 and 12) the first is included in the training data and the second one in the application data. 
\begin{table} \centering
  \caption{Recording times and statistics on the passing and entering traffic for all data sets. Data was recorded in 2019 on May 13 (day 1), May 14 (day 2), May 15 (day 3) and June 21 (day 4). Time in the table refers to the start of the recordings. The number and the average speed for both passing and entering vehicles are shown.}\label{tab:datasets}
  \begin{tabular}{r r r l r r r r r r}\toprule
    set & \multicolumn{2}{c}{time} & period & \multicolumn{2}{c}{passing vehicles} & \multicolumn{2}{c}{entering vehicles} & \multicolumn{2}{c}{coupling delay}\\
    & & & & n & av. speed& n& av. speed & $\tau^2$ & $\tau^3$ \\\midrule
    1 & day 1& 5:10:30 & 5 min, 25 sec & 281 & 74.9 km/h & 52 & 66.1 km/h & 0.75 & 7.00 \\ 
    2 & day 1& 5:15:57 & 5 min, 26 sec & 301 & 74.3 km/h & 71 & 63.5 km/h & -0.25 & 6.00 \\ 
    3 & day 1& 5:21:24 & 5 min, 32 sec & 298 & 71.0 km/h & 71 & 62.3 km/h & -0.25 & 6.75 \\ 
    4 & day 1& 5:43:57 & 5 min, 26 sec & 312 & 67.3 km/h & 73 & 56.3 km/h & 0.50 & 9.75\\ 
    5 & day 1& 5:49:25 & 5 min, 26 sec & 288 & 65.0 km/h & 76 & 53.3 km/h & 0.00 & 9.25\\ 
    6 & day 1& 5:54:52 & 5 min, 27 sec & 286 & 67.3 km/h & 68 & 56.2 km/h & 0.00 & 8.75\\ 
    7 & day 1& 6:00:22 & 1 min, 38 sec & 28 & 68.2 km/h & 3 & 49.4 km/h &  4.50 & 13.50\\ 
    8 & day 1& 6:04:34 & 5 min, 26 sec & 282 & 67.7 km/h & 89 & 56.3 km/h & 1.25 & 9.50\\ 
    9 & day 1& 6:38:58 & 5 min, 26 sec & 269 & 67.4 km/h & 64 & 56.9 km/h & 0.00 & 9.00\\ 
    10 & day 1& 6:44:26 & 1 min, 25 sec & 63 & 48.8 km/h & 11 & 32.2 km/h & -5.00 & 9.75\\ 
    11 & day 2& 15:10:37 & 2 min, 45 sec & 85 & 74.4 km/h & 18 & 61.7 km/h & -0.50 & 7.25\\ 
    12 & day 2& 15:26:35 & 5 min, 27 sec & 215 & 39.2 km/h & 33 & 33.3 km/h &  5.00 & 17.25\\ 
    13 & day 2& 15:44:39 & 5 min, 27 sec & 233 & 66.3 km/h & 50 & 54.7 km/h & 0.00 & 9.00\\ 
    14 & day 2& 15:50:07 & 5 min, 27 sec & 235 & 73.0 km/h & 54 & 61.0 km/h & 0.00 & 8.50\\ 
    15 & day 2& 15:55:34 & 2 min, 22 sec & 107 & 73.3 km/h & 29 & 63.8 km/h & -0.25 & 8.25\\ 
    16 & day 2& 16:02:44 & 5 min, 27 sec & 250 & 73.5 km/h & 69 & 62.6 km/h & -0.25 & 7.50\\ 
    17 & day 2& 16:08:13 & 5 min, 26 sec & 244 & 74.2 km/h & 59 & 61.4 km/h & -0.25 & 7.75\\ 
    18 & day 2& 16:19:24 & 3 min, 6 sec & 149 & 67.3 km/h & 37 & 61.4 km/h & -0.25 & 7.75\\ 
    19 & day 3& 6:07:22 & 5 min, 27 sec & 300 & 65.2 km/h & 89 & 54.4 km/h & -0.75 & 9.00\\ 
    20 & day 3& 6:12:51 & 5 min, 26 sec & 278 & 63.8 km/h & 83 & 52.2 km/h & 4.50 & 13.25\\ 
    21 & day 3& 7:26:13 & 2 min, 39 sec & 66 & 67.4 km/h & 7 & 58.7 km/h & 0.75 & 9.75\\ 
    22 & day 4& 15:28:14 & 5 min, 26 sec & 221 & 72.5 km/h & 27 & 60.5 km/h & 0.00 & 8.50\\ 
    23 & day 4& 15:33:42 & 5 min, 26 sec & 247 & 72.3 km/h & 20 & 61.9 km/h & 0.00 & 8.50\\ 
    24 & day 4& 15:39:09 & 5 min, 26 sec & 230 & 73.7 km/h & 15 & 62.5 km/h & 0.00 & 8.25\\ 
    25 & day 4& 15:44:36 & 2 min, 48 sec & 76 & 67.8 km/h & 4 & 58.9 km/h & -4.50 & 4.75\\ 
    26 & day 4& 16:06:02 & 2 min, 59 sec & 79 & 68.6 km/h & 9 & 63.231687 & 0.00 & 8.50\\ 
    27 & day 4& 16:11:53 & 5 min, 26 sec & 231 & 77.3 km/h & 31 & 62.3 km/h & -5.00 & 2.75\\ 
    28 & day 4& 16:33:19 & 5 min, 27 sec & 231 & 74.3 km/h & 29 & 62.2 km/h & 0.00 & 7.75\\ 
    29 & day 4& 16:38:47 & 5 min, 27 sec & 244 & 68.2 km/h & 38 & 57.1 km/h & 0.00 & 8.00\\ 
    30 & day 4& 16:44:15 & 5 min, 27 sec & 248 & 73.0 km/h & 34 & 62.0 km/h & -3.00 & 5.00\\ 
    31 & day 4& 16:49:45 & 4 min, 34 sec & 171 & 68.6 km/h & 14 & 54.7 km/h & 0.00 & 8.50\\
    \bottomrule
  \end{tabular}
\end{table}
\bibliographystyle{abbrvurl} 
\bibliography{traffic.bib}
\end{document}